\documentclass[10pt]{amsart} 
\title{A Markov chain on permutations which projects to the PASEP}
\author{Sylvie Corteel and Lauren K. Williams}

\address{CNRS LRI, Universit\'e Paris-Sud, B\^atiment 490, 91405 Orsay
Cedex France} 
\email{Sylvie.Corteel@lri.fr}
\address{Harvard University, Cambridge}
\email{lauren@math.harvard.edu}

\subjclass[2000]{Primary 05E10; Secondary 82B23, 60C05}

\keywords{}

\usepackage{pstricks, pst-node, amssymb}

\psset{unit=1pt, arrowsize=4pt, linewidth=.7pt}
\psset{linecolor=blue}
\newgray{grayish}{.90}
\newrgbcolor{embgreen}{0 .5 0}
\def\Le{\hbox{\rotatedown{$\Gamma$}}}
\def\vblack(#1, #2)#3{\cnode*[linecolor=black](#1, #2){3}{#3}}
\def\vwhite(#1,#2)#3{\cnode[linecolor=black,fillcolor=white,fillstyle=solid](#1,
#2){3}{#3}}
\countdef\x=23
\countdef\y=24
\countdef\z=25
\countdef\t=26

\def\tbox(#1,#2)#3{
\x=#1 \y=#2
\multiply\x by 12
\multiply\y by 12
\z=\x \t=\y
\advance\z by 12
\advance\t by 12
\psline(\x,\y)(\x,\t)(\z,\t)(\z,\y)(\x,\y)
\advance\x by 6
\advance\y by 6
\rput(\x,\y){{\bf #3}}}

\usepackage{amsmath, amsthm, amssymb, amsbsy}
\usepackage{amsfonts, latexsym, stmaryrd, amscd, xy}
\usepackage[mathscr]{eucal}
\usepackage{epsfig}

\newtheorem{theorem}{Theorem}[section]

\newtheorem{proposition}[theorem]{Proposition}
\newtheorem{lemma}[theorem]{Lemma}

\newtheorem{example}[theorem]{Example}
\newtheorem{corollary}[theorem]{Corollary}

\newtheorem{remark}[theorem]{Remark}

\newtheorem{definition}[theorem]{Definition}

\usepackage{amsfonts}
\usepackage{xy}
\usepackage{amssymb}
\usepackage{amsmath}
\newcommand{\T}{\mathcal{T}}
\newcommand{\SSS}{\mathcal{S}}

\newcommand{\ttt}{\tau}

\newcommand{\Q}{\mathcal{Q}}
\newcommand{\Sset}{\mathcal{S}}

\DeclareMathOperator{\half-perimeter}{half-perim}
\DeclareMathOperator{\Norm}{Norm}
\DeclareMathOperator{\wt}{wt}
\DeclareMathOperator{\rk}{rk}

\DeclareMathOperator{\prob}{prob}
\DeclareMathOperator{\pr}{pr}

\newcommand{\thmrefer}[1]{\renewcommand\thetheorem
  {\protect\ref{#1}}\addtocounter{theorem}{-1}}

\xyoption{all}
\CompileMatrices

\begin{document}


\begin{abstract}

The partially asymmetric exclusion process (PASEP) is an important model
from statistical mechanics which describes a system of interacting 
particles hopping left and right on a one-dimensional lattice of $N$
sites.  It is partially asymmetric in the sense that the probability 
of hopping left is $q$ times the probability of hopping right.  
Additionally, particles may enter from the left with probability 
$\alpha$ and exit to the right with probability $\beta$.

It has been observed that the (unique) stationary distribution of the PASEP
has remarkable connections to combinatorics -- see for example the papers of 
Derrida {\it et al} \cite{Derrida0, Derrida1}, Duchi and Schaeffer \cite{jumping},
Corteel \cite{Corteel}, and Shapiro and Zeilberger \cite{SZ}. 
Most recently we proved \cite{CW} that in fact
the (normalized) probability of being in a particular state of the PASEP can be viewed as 
a certain weight generating function for {\it permutation tableaux} of a fixed 
shape.  (This result implies the previous combinatorial results.)  
However, our proof relied on the {\it matrix ansatz} of 
Derrida {\it et al} \cite{Derrida1}, and hence did not give an intuitive
explanation of why one should expect the steady state distribution of the PASEP to 
involve such nice combinatorics.

In this paper we 
define a Markov chain -- 
which we call the PT chain -- on the set of permutation 
tableaux which {\it projects} to the PASEP, in a sense which we shall
make precise. 
This gives a new proof of the main 
result of \cite{CW} which bypasses the matrix ansatz altogether.  
Furthermore, via the bijection of \cite{SW}, the PT 
chain can also be viewed as a Markov chain on the symmetric group.
Another nice feature of the PT chain is that it possesses a certain symmetry
which extends the {\it particle-hole symmetry} of the PASEP.
\end{abstract}

\maketitle

\section{Introduction}

The partially asymmetric exclusion process (PASEP) 
is an important model from statistical mechanics which is quite  simple 
but surprisingly rich: it exhibits 
boundary-induced phase transitions, spontaneous symmetry 
breaking, and phase separation.  The PASEP is regarded as a primitive
model for biopolymerization \cite{bio}, traffic flow \cite{traffic}, 
and formation of shocks \cite{shock};
it also appears in a kind of sequence alignment problem in 
computation biology \cite{comput}.

In brief, the PASEP  
describes a system of particles hopping left and right on a one-dimensional lattice of $N$ sites.  Particles may enter the system from 
the left with a rate $\alpha dt$ and may exit the system to the right 
at a rate $\beta dt$.  The probability of hopping left is $q$ times the 
probability of hopping right.

It has been observed that the (unique) stationary distribution of the PASEP
has remarkable connections to combinatorics.  
Derrida {\it et al} \cite{Derrida0, Derrida1} 
proved a connection to Catalan and Narayana numbers in the 
case where $q=0$ (TASEP) and $\alpha=\beta=1$;
Duchi and Schaeffer \cite{jumping} gave a combinatorial 
explanation of this result by 
constructing a new Markov chain on ``complete configurations"
(enumerated by Catalan numbers) that projects to the TASEP, for 
general $\alpha$ and $\beta$.
Subsequently Corteel \cite{Corteel} proved a connection of the PASEP
to the $q$-Eulerian numbers of 
\cite{Williams}, thereby generalizing  
the result of Derrida {\it et al} to include the case where again 
$\alpha=\beta=1$ but $q$ is general. 
Shortly thereafter we proved \cite{CW} a 
much stronger result in the case of general 
$\alpha, \beta$ and $q$, showing that in fact 
the (normalized) probability of being in a particular state of the PASEP can be viewed as 
a certain {\it weight generating function} for {\it permutation tableaux} 
(certain $0-1$ tableaux) of a fixed 
shape -- this is a Laurent polynomial in $\alpha$, $\beta$, and $q$.  
However, our proof relied on the {\it matrix ansatz} of 
Derrida {\it et al} \cite{Derrida1}, and hence did not give an intuitive
explanation of why one should expect the steady state of the PASEP to 
be related to such nice combinatorics.

The goal of this paper is to construct a Markov chain on permutation tableaux
(which we call the PT chain)
which {\it projects} to the PASEP: that is,
after projection via a certain surjective map between the spaces of states,
a walk on the state diagram of the PT chain is indistinguishable from 
a walk on the state diagram of the PASEP.  
The steady state distribution of the PT chain
has the nice property that the (normalized) probability 
of being in a particular state (i.e.\ a permutation tableau) is the 
{\it weight} of that permutation tableau -- this is a Laurent {\it monomial} in 
$\alpha$, $\beta$, and $q$.
Our construction generalizes the work of 
Duchi and Schaeffer \cite{jumping} (whose work can be viewed as the $q=0$ case of 
ours), and gives a new proof of the main result of \cite{CW}.
Additionally by using the bijection of \cite{SW}, we can view the PT chain 
as a Markov chain on the symmetric group.  
Finally, the PT chain  possesses a certain
symmetry which extends the {\it particle-hole symmetry} of the PASEP:
this is a graph-automorphism on the state diagram of the PT chain which
is an involution.

The structure of this paper is as follows.  
In Section \ref{setup} we define the PASEP.
In Section \ref{PermTableaux}
we define permutation tableaux, certain $0-1$ tableaux which 
are naturally in bijection
with permutations.  Section \ref{Markov-chain} 
defines the PT chain, introduces the notion of {\it projection} of Markov
chains, and proves a tight relationship between  the steady state
distributions of two Markov chains when one projects to the other.  
It will be clear from the definition that 
the PT chain projects to the PASEP.
Section \ref{Proof}
states and proves our main result about the steady 
state distribution of the PT 
chain, thus giving a new combinatorial proof of our main result from \cite{CW}.
Section \ref{Perm-Section} recalls the definition 
$\Phi$ from \cite{SW}, and uses it to describe both the PT chain 
as a Markov chain on permutations.
Finally, Section \ref{involution} describes an involution on the 
state-diagram of the PT chain which extends the 
particle-hole symmetry.

It would be interesting to explore whether the PT chain has any
physical significance, and 
whether this larger chain may shed some insight on 
the PASEP itself.

\textsc{Acknowledgments:}  We are grateful to Mireille 
Bousquet-Melou, Bernard Derrida, Persi Diaconis, Philippe Duchon, 
Daniel Ford, 
and Xavier Viennot 
for their comments and encouragement.  Additionally, we thank
the referee for advice that enabled us to considerably simplify 
our arguments concerning the involution.

\section{The PASEP}\label{setup}

In the physics literature, the PASEP is defined as follows.
\begin{definition}
We are given a one-dimensional
lattice of $N$ sites, such that each site $i$ ($1 \leq i \leq N)$
is either occupied by a particle ($\tau_i=1$) or is empty
($\tau_i=0$).  At most one particle may occupy a given site.
During each infinitesimal time interval $dt$, each particle 
in the system has a probability $dt$ of jumping to the next 
site on its right (for particles on sites $1 \leq i \leq N-1$) and 
a probability $q dt$ of jumping to the next site on its left
(for particles on sites $2 \leq i \leq N$).  Furthermore, a particle
is added at site $i=1$ with probability $\alpha dt$ if site $1$
is empty and a particle is removed from site $N$ with probability 
$\beta dt$ if this site is occupied.  
\end{definition}

\begin{remark}
Note that we will sometimes denote a state of the PASEP as a 
word in $\{0,1\}^N$ and sometimes as a word in $\{\circ,\bullet\}^N$.
In the latter notation, the symbol $\circ$ denotes the absence of 
a particle, which one can also think of as a white particle.
\end{remark}

It is not too hard to see \cite{jumping} that our previous formulation
of the PASEP is equivalent to the following  
discrete-time Markov chain.

\begin{definition}
Let $B_N$ be the set of all $2^N$ words in the
language $\{\circ, \bullet\}^*$.
The PASEP is the Markov chain on $B_N$ with
transition probabilities:
\begin{itemize}                                                                
\item  If $X = A\bullet \circ B$ and
$Y = A \circ \bullet B$ then
$P_{X,Y} = \frac{1}{N+1}$ (particle hops right) and
$P_{Y,X} = \frac{q}{N+1}$ (particle hops left).
\item  If $X = \circ B$ and $Y = \bullet B$
then $P_{X,Y} = \frac{\alpha}{N+1}$ (particle enters from left).
\item  If $X = B \bullet$ and $Y = B \circ$
then $P_{X,Y} = \frac{\beta}{N+1}$ (particle exits to the right).
\item  Otherwise $P_{X,Y} = 0$ for $Y \neq X$
and $P_{X,X} = 1 - \sum_{X \neq Y} P_{X,Y}$.
\end{itemize}
\end{definition}

\begin{figure}[h]
\centering
\includegraphics[height=1.3in]{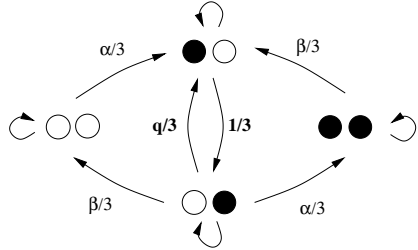}
\caption{The state diagram  of the PASEP for $N=2$}
\label{states}
\end{figure}

See Figure \ref{states} for an
illustration of the four states, with transition probabilities,
for the case $N=2$.

In the long time limit, the system reaches a steady state where all 
the probabilities $P_N(\ttt_1, \ttt_2, \dots , \ttt_N)$ of finding
the system in configurations $(\ttt_1, \ttt_2, \dots , \ttt_N)$ are
stationary, i.e.\ satisfy 
\begin{equation*}
\frac{d}{dt} P_N(\tau_1, \dots , \ttt_N) = 0.
\end{equation*}
Moreover, the stationary distribution is unique \cite{Derrida1}, as shown
by Derrida {\it et al}.

The question is now to solve for the probabilities 
$P_N(\ttt_1, \dots , \ttt_N)$.  For convenience, we 
work with unnormalized weights $g_N(\ttt_1, \dots , \ttt_N)$ which are
equal to the $P_N(\ttt_1, \dots , \ttt_N)$ up to a constant:
\begin{equation*}
P_N(\ttt_1, \dots , \ttt_N) = g_N(\ttt_1, \dots , \ttt_N)/Z_N,
\end{equation*}
where 
$Z_N$ is the {\it partition function}
$\sum_{\tau} g_N(\ttt_1, \dots , \ttt_N)$.  The sum defining
$Z_N$ is 
over all possible configurations $\ttt \in \{0,1\}^N$.

\begin{remark}\label{particle-hole}
There is an obvious {\it particle-hole symmetry} \cite{Derrida1} in the PASEP:
since (black) particles enter at the left with probability 
$\alpha$ and exit to the right with probability $\beta$, it is equivalent
to saying that holes (or white particles) are injected at the right with probability
$\beta$ and are removed at the left end with probability $\alpha$.  

Let us define $\overline{(\tau_1, \dots , \tau_N)} = (1-\tau_N, 1-\tau_{N-1}, \dots , 1-\tau_1)$.
Clearly this operation is an involution on states of the PASEP.
Because  
of the particle-hole symmetry, one always has 
that $g_N^{q,\alpha,\beta}(\tau) = 
g_N^{q,\beta,\alpha}(\overline{\tau})$.

\end{remark}

\section{Connection with permutation tableaux}\label{PermTableaux}

We define a {\em partition} $\lambda = (\lambda_1, \dots,
\lambda_K)$ to be a weakly decreasing sequence of {\it nonnegative}
integers. For a partition $\lambda$, where $\sum \lambda_i = m$, the
{\em Young diagram} $Y_\lambda$ of shape $\lambda$ is a left-justified
diagram of $m$ boxes, with $\lambda_i$ boxes in the $i$th row.
We define the {\it half-perimeter} of $\lambda$ or $Y_\lambda$ to be the  
the number of rows plus the number of columns.  
The {\it length} of a row or column of a Young diagram is the number
of boxes in that row or column.  Note that 
we will allow a row to have length $0$.

\begin{figure}[h]
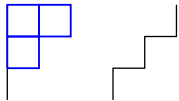

\pspicture(0,0)(40,40)        
\psline[linecolor=black,linewidth=0.5pt]{-}(0,0)(0,36)
\tbox(0,2){}                 
\tbox(0,1){}                 
\tbox(1,2){}                 
\endpspicture
\pspicture(0,0)(40,40)
\psline[linecolor=black,linewidth=0.5pt]{-}(0,0)(0,12)(12,12)(12,24)(24,24)(24,36)
\endpspicture
\caption{Young diagram and path for $\lambda=(2,1,0)$}
\label{path}
\end{figure}

We will often identify
a Young diagram $Y_\lambda$ of half-perimeter $t$ 
with the lattice path $p(\lambda)$ of length $t$
which takes unit steps south and west,
beginning at
the north-east corner of $Y_\lambda$ and ending at
the south-west corner.
Note that such a lattice path always begins with a step south.
Figure \ref{path} shows the path corresponding to the 
Young diagram of shape $(2,1,0)$.

If $\tau \in \{0,1\}^{N}$, we define a Young diagram 
$\lambda(\tau)$ of half-perimeter $N+1$ as follows.  First we define 
a path $p = (p_1, \dots , p_{N+1}) \in \{S,W\}^{N+1}$ such that $p_1=S$, and 
$p_{i+1}=S$ if and only if $\tau_i = 1$.  We then define
$\lambda(\tau)$ to be the partition associated to this path $p$.
This map is clearly a bijection between the set of Young diagrams
of half-perimeter $N+1$ and the set of $N$-tuples in $\{0,1\}^N$,
and we denote the inverse map similarly: 
given a Young diagram $\lambda$ of half-perimeter $N+1$, we define $\tau(\lambda)$ to be
the corresponding $N$-tuple.

As in \cite{SW}, we
define a {\em permutation tableau} $\T$ to be a partition $\lambda$
together with a filling of the boxes of $Y_\lambda$ with $0$'s and
$1$'s such that the following properties hold:
\begin{enumerate}
\item Each column of the rectangle contains at least one $1$.
\item There is no $0$ which has a $1$ above it in the same column
{\em and} a $1$ to its left in the same row.
\end{enumerate}

We call such a filling a {\em valid} filling of $Y_\lambda$.  

\begin{remark}
Permutation tableaux are closely 
connected to total positivity for the Grassmannian
\cite{Postnikov, Williams}.   
More precisely, 
if we forget the requirement (1) above we recover the definition of
a $\Le$-diagram, 
an object which represents
a cell in the totally nonnegative part of the Grassmannian.
It would be interesting to explore
whether there is a connection between total positivity and the PASEP.
\end{remark}

\begin{remark}
Sometimes we will depict permutation tableaux slightly differently,
replacing the $1$'s with black dots and omitting the $0$'s entirely,
as in Figure \ref{Markov-states}.
\end{remark}

Note that the second requirement above can be rephrased in the following
way.  Read the columns of a permutation tableau $\T$ from right to left.  
If in any column we have a $0$ which lies beneath some $1$, then all 
entries to the left of $0$ (which are in the same row) must also
be $0$'s.

We  now define a few statistics on permutation tableaux.
A $1$ in a tableau $\T$ is {\it topmost} if it has only $0$'s above it;
a $1$ in $\T$ is {\it superfluous} if it is not topmost; and a $1$ is 
{\it necessary} if it is the unique $1$ in its column.
We define the {\it rank} $\rk(\T)$ of a permutation tableau $\T$ 
to be the number of superfluous $1$'s.  (Therefore $\rk(\T)$ is equal
to the total number of $1$'s in the filling minus the number of columns.)
We define $f(\T)$ to be the number of $1$'s in the first
row of $\T$.
We say that a zero in a permutation tableau is 
{\it restricted} if 
there is a one above it in the same column. 
And we say that a row is 
{\it unrestricted}
if it does not contain a restricted entry.  
Define $u(\T)$ to be the number of unrestricted rows of $\T$ 
minus $1$.  (We subtract $1$ since the top row of a tableau
is always unrestricted.)

\begin{figure}[h]
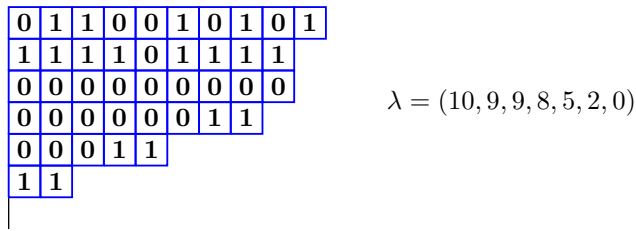

\pspicture(-95,-12)(230,94)
\rput(190,36)
{$\begin{array}{l}
\lambda=(10,9,9,8,5,2,0)
\end{array}$}
\psline[linecolor=black,linewidth=0.5pt]{-}(0,-12)(0,0)
\tbox(0,0){1}
\tbox(1,0){1}
\tbox(0,1){0}
\tbox(1,1){0}
\tbox(2,1){0}
\tbox(3,1){1}
\tbox(4,1){1}
\tbox(0,2){0}
\tbox(1,2){0}
\tbox(2,2){0}
\tbox(3,2){0}
\tbox(4,2){0}
\tbox(5,2){0}
\tbox(6,2){1}
\tbox(7,2){1}
\tbox(0,3){0}
\tbox(1,3){0}
\tbox(2,3){0}
\tbox(3,3){0}
\tbox(4,3){0}
\tbox(5,3){0}
\tbox(6,3){0}
\tbox(7,3){0}
\tbox(8,3){0}
\tbox(0,4){1}
\tbox(1,4){1}
\tbox(2,4){1}
\tbox(3,4){1}
\tbox(4,4){0}
\tbox(5,4){1}
\tbox(6,4){1}
\tbox(7,4){1}
\tbox(8,4){1}
\tbox(0,5){0}
\tbox(1,5){1}
\tbox(2,5){1}
\tbox(3,5){0}
\tbox(4,5){0}
\tbox(5,5){1}
\tbox(6,5){0}
\tbox(7,5){1}
\tbox(8,5){0}
\tbox(9,5){1}
\endpspicture
\caption{A permutation tableau}
\label{PermTab}
\end{figure}

Figure \ref{PermTab} gives an example of a permutation tableau $\T$
with rank $19-10=9$ and half-perimeter $17$, such that 
$u(\T) = 3$ and $f(\T) = 5$.

We
define the {\it weight}
of a tableau $\T$ to be the monomial 
$\wt(\T):=q^{\rk(\T)} \alpha^{-f(\T)} \beta^{-u(\T)}$,
and we define 
$F_{\lambda}(q)$ to be the (Laurent) polynomial $\sum_{\T} \wt(\T)$, 
where the 
sum ranges over all permutation tableaux $\T$ of shape
$\lambda$.  

Our main result of  \cite{CW} was the following.

\begin{theorem}\label{oldtheorem}
Fix $\ttt = (\ttt_1, \dots , \ttt_N) \in \{0,1\}^N$, and 
let $\lambda := \lambda(\tau)$.  (Note that $\half-perimeter (\lambda)=N+1$.)
The probability of finding the PASEP in configuration 
$(\ttt_1, \dots , \ttt_N)$ in the steady state is 
\begin{equation*}
\frac{F_{\lambda}(q)}{Z_N}.
\end{equation*}
Here, $F_{\lambda}(q)$ is the weight-generating function for 
permutation tableaux of shape $\lambda$. 
Moreover, the partition function $Z_N$ for the PASEP 
is equal to the weight-generating function for all permutation
tableaux of half-perimeter $N+1$. 
\end{theorem}

However, our proof relied on the {\it matrix ansatz} of 
Derrida {\it et al}.  In this paper we 
will give another proof of that result, which 
bypasses the matrix ansatz and gives a better explanation of the 
connection between permutation tableaux and the PASEP.  


\section{The PT 
chain}\label{Markov-chain}

In this section we present our main construction, 
a Markov chain
on permutation tableaux which projects to the PASEP.
We will call this chain the 
{\it PT chain}. 

\begin{definition}
We define a projection operator $\pr$ which projects
a state of the PT, i.e.\ a permutation tableau, to a state
of the PASEP.  If $\T$ is a permutation tableau of shape $\lambda$
and half-perimeter  
$N+1$, then we define 
$\pr(\T):= \tau(\lambda)$.  This is a state of the PASEP with $N$ sites.
\end{definition}

\begin{figure}[h]
\centering
\includegraphics[height=3.5in]{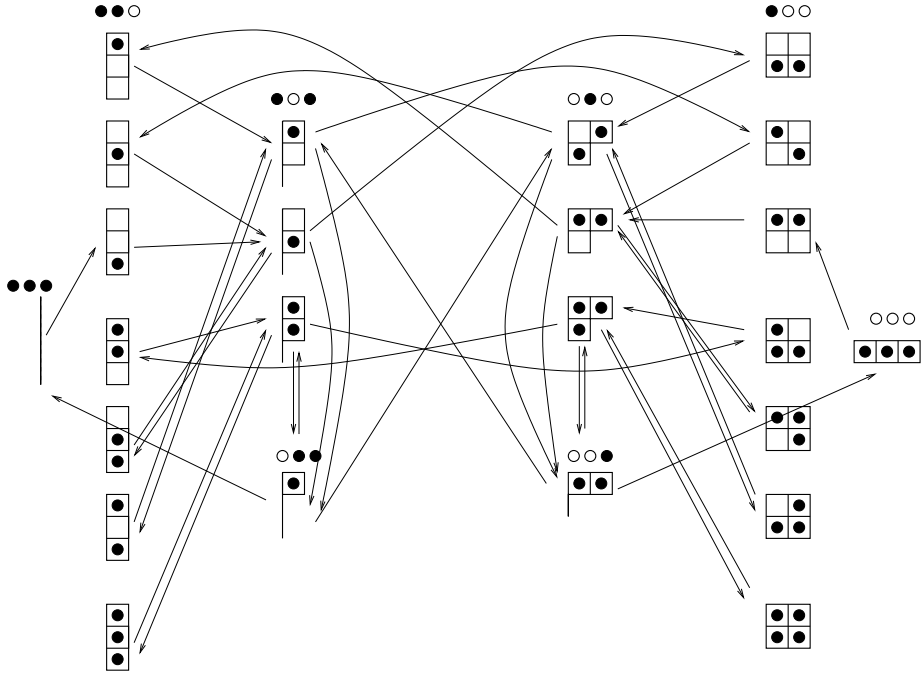}
\caption{The state diagram of the PT chain for $N=3$}
\label{Markov-states}
\end{figure}

Before defining the PT chain, we show an example: the 
state diagram of the chain 
for $N=3$.  In Figure \ref{Markov-states}, the $24=4!$ states
of the PT chain are arranged into $8=2^3$ groups according 
to partition shape (four of size $1$,
two of size $3$, two of size $7$).  All elements of a 
fixed  group of tableaux project to the same state of the PASEP,
depicted just above that group.  We have 
not included the transition probabilities in Figure 
\ref{Markov-states}, but they are defined in the obvious manner:
if there is a transition $\SSS\to T$ in the PT chain,
then
$\prob_{PT}(\SSS \to T) = \prob_{PASEP}(\pr(\SSS) \to \pr(T))$.
Finally, observe that there is a reflective left-right symmetry
in the figure.

We now define all possible transitions in the PT chain.
There are four 
kinds of transitions, which correspond to the four kinds 
of transitions in the PASEP.  The transition probabilities
in the PT chain will be defined in accordance with the 
transition probabilities in the PASEP:
if there is a transition $\SSS\to T$ in the PT chain,
then
$\prob_{PT}(\SSS \to T) := \prob_{PASEP}(\pr(\SSS) \to \pr(T))$.

In what follows,
we will assume that $\SSS$ is a permutation tableau of half-perimeter $N+1$,
whose shape is $\lambda = (\lambda_1, \dots , \lambda_m, \dots , 
\lambda_t)$ 
where $\lambda_1\geq \dots \geq \lambda_m>0$, and $\lambda_r=0$ for $r>m$.

\subsection{Particle enters from the left}\label{enterleft}

If the rightmost column of $\SSS$ has length $1$, then
there is a transition in the PT chain 
from $\SSS$ that corresponds to a particle 
entering from the left in the PASEP. 

\begin{figure}[h]
\centering
\includegraphics[height=1.1in]{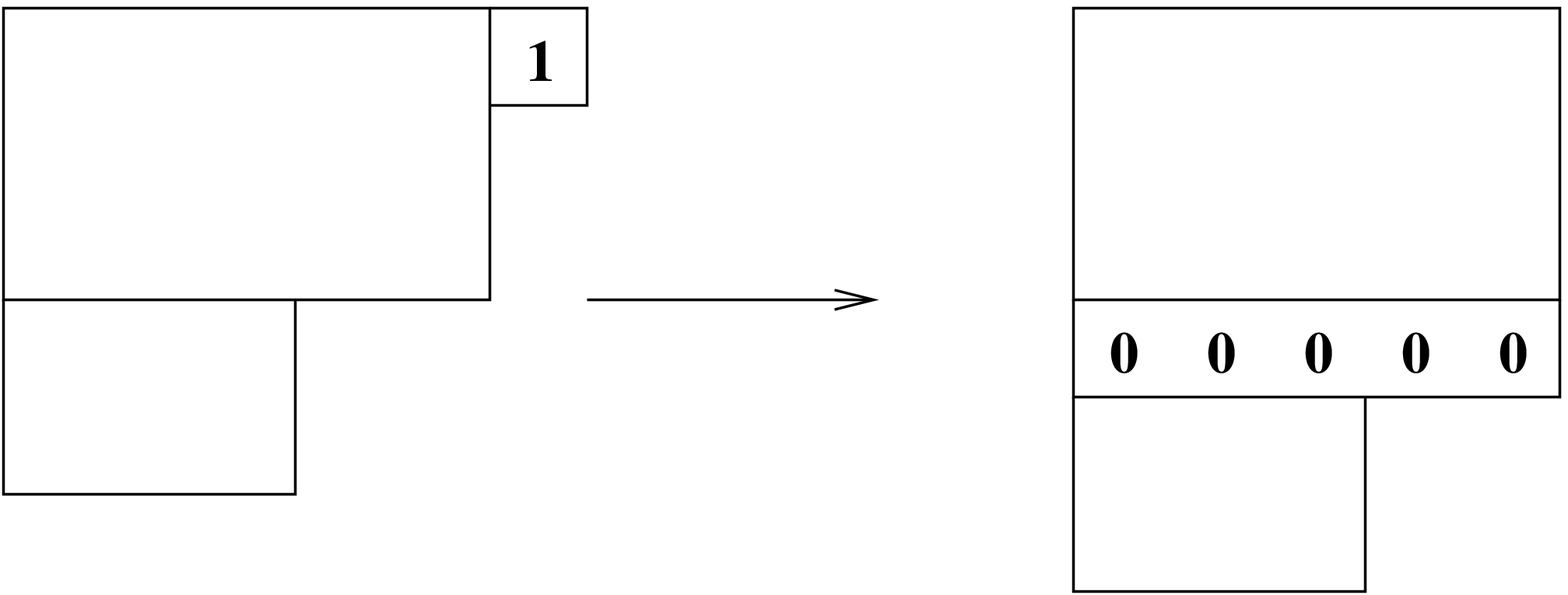}
\caption{}
\label{enter}
\end{figure}

We now define a new permutation tableau $\T$ as follows:
delete the rightmost column of $\SSS$ and add a new all-zero 
row of length $\lambda_1 -1$ to $\SSS$, inserting it as far south 
as possible (subject to the constraint that the lengths of the rows
of a permutation tableau must weakly decrease). Clearly 
adding an all-zero row in this way results in a new permutation
tableau, since there is no way to introduce the forbidden pattern
(condition (2) in the definition of permutation tableau),
and each column will still contain at least one $1$.
See Figure 
\ref{enter}.

We define $\prob(\SSS \to \T) = \frac{\alpha}{N+1}$.  Since we have
removed a $1$ in the top row but have not affected the unrestricted 
rows or superfluous $1$'s, we have that
$\wt(\T) = \alpha \cdot \wt(\SSS)$, and therefore 
$\wt(\SSS) \cdot \prob(\SSS \to \T) = \frac{\wt(\T)}{N+1}$.

\subsection{Particle hops right}\label{hopright}

If some row $\lambda_j > \lambda_{j+1}$ in $\SSS$ (resp. if $\lambda_t > 0$),
then there 
is a transition in PT from $\SSS$ that corresponds to the
$(j-1)$st black particle (resp. $(t-1)$st black particle) in $\pr(\SSS)$ 
hopping to the right in the PASEP.

We now define a new permutation tableau $\T$ as follows, 
based on the rightmost entry of the $j$th row of $\SSS$.

\subsubsection{Case 1}
\begin{figure}[h]
\centering
\includegraphics[height=1.1in]{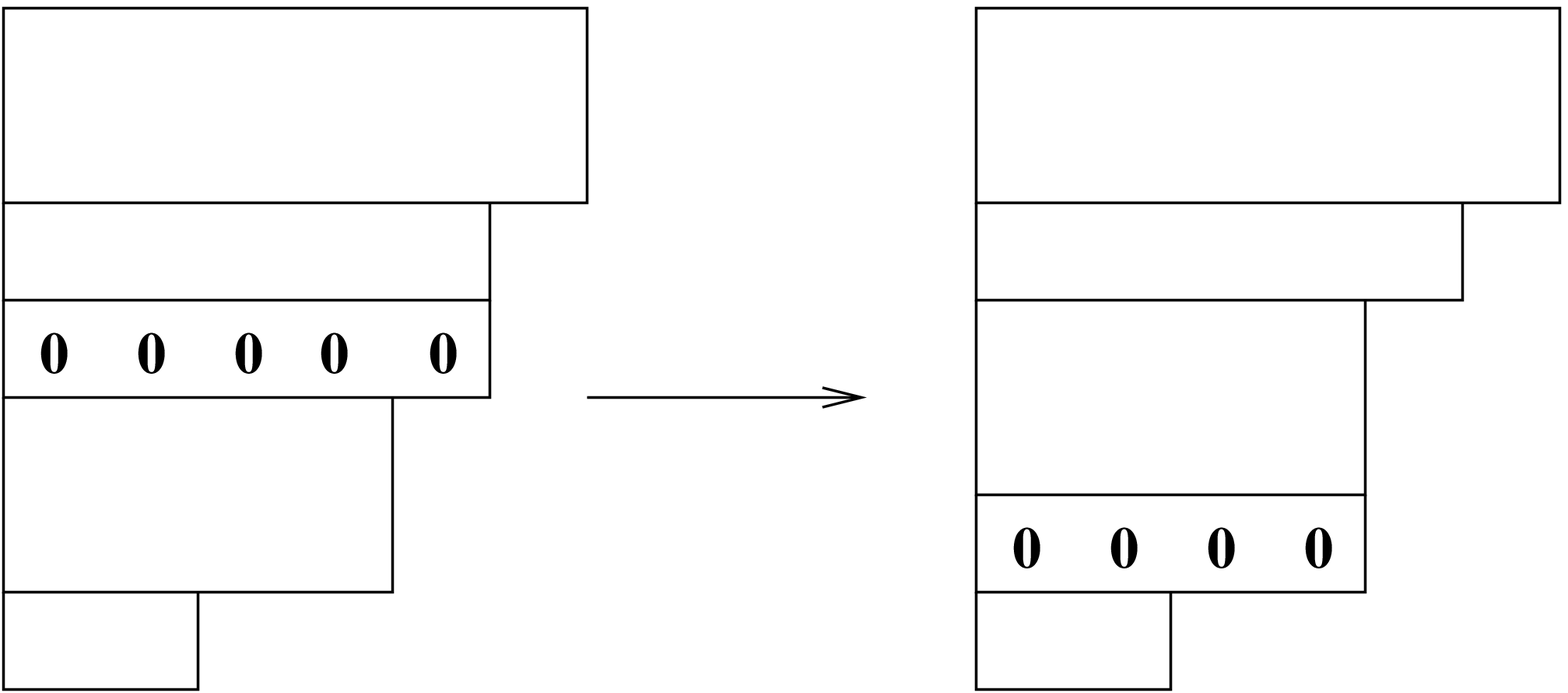}
\caption{}
\label{Right1}
\end{figure}

Suppose that the rightmost entry of the $j$th row is a $0$.
This forces the $j$th row to contain only $0$'s.
Then we define a new tableau $\T$ by deleting the $j$th row of $\SSS$
and adding a new row of $\lambda_j - 1$ $0$'s, inserting it as far
south as possible.  Clearly adding and deleting all-zero rows
in this fashion results in a permutation tableau, so this operation
is well-defined.
See Figure \ref{Right1}.

We define $\prob(\SSS \to \T) = \frac{1}{N+1}$.  
If $\lambda_j > 1$ then we have not affected the number
of superfluous $1$'s, unrestricted rows, or $1$'s in the top row, and 
so
$\wt(\T) = \wt(\SSS)$.  It follows that 
$\wt(\SSS) \cdot \prob(\SSS \to \T) = \frac{\wt(\T)}{N+1}$.

In the special case that  $\lambda_j = 1$, then 
$\wt(\T) = {\beta}^{-1} \wt(\SSS)$, and therefore 
$\wt(\SSS) \cdot \prob(\SSS \to \T) = \frac{\beta \wt(\T)}{N+1}$.
Note that if this special case occurs then $\pr(\T)$ ends with 
a black particle.  And given such a $\T$, there is only one 
such $\SSS$ with such a transition $\SSS \to \T$.

\subsubsection{Case 2}
\begin{figure}[h]
\centering
\includegraphics[height=1.1in]{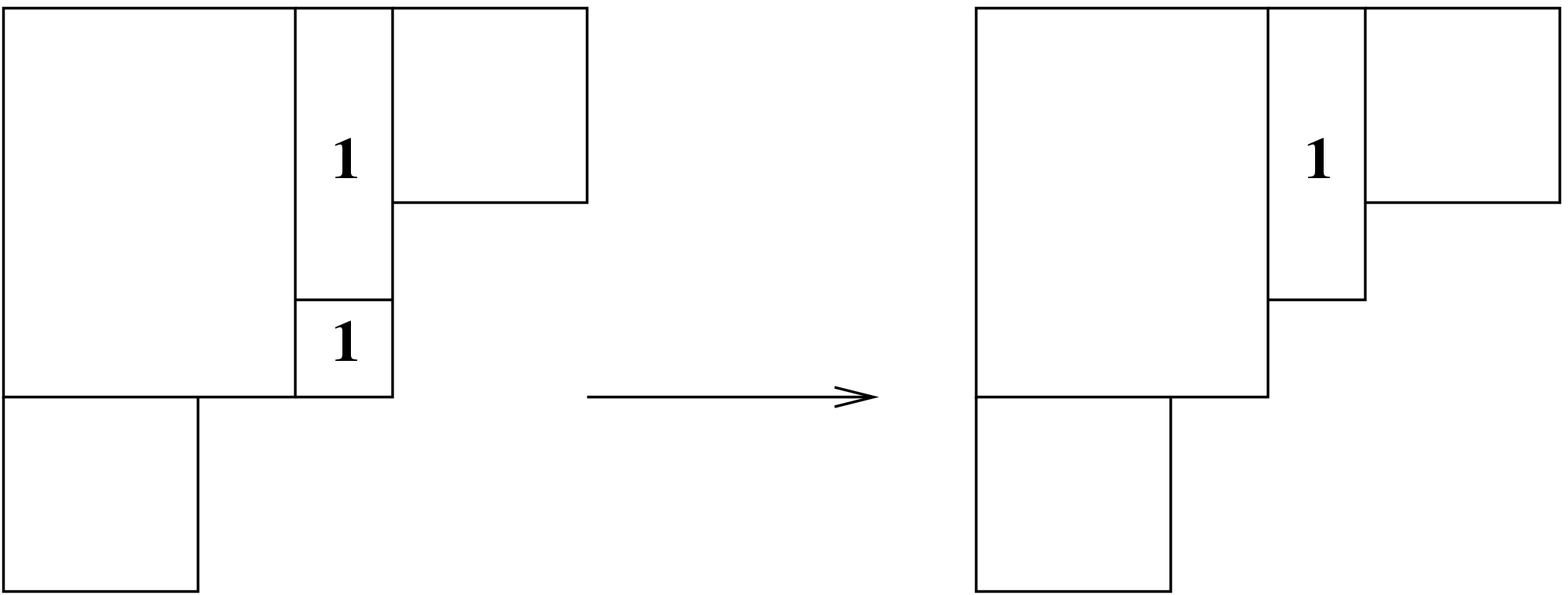}
\caption{}
\label{Right2}
\end{figure}

Suppose that the rightmost entry of the $j$th row of $\SSS$ is a 
superfluous $1$.
Then we define a new tableau $\T$ by deleting that $1$ from $\SSS$;
because the column containing that $1$ had at least two $1$'s
to begin with, our new tableau will be a permutation tableau.
See Figure \ref{Right2}.

We define $\prob(\SSS \to \T) = \frac{1}{N+1}$.  
Our operation has affected only the number of superfluous $1$'s,
so
$\wt(\T) = q^{-1} \wt(\SSS)$, and therefore 
$\wt(\SSS) \cdot \prob(\SSS \to \T) = \frac{q \wt(\T)}{N+1}$.

\subsubsection{Case 3}
\begin{figure}[h]
\centering
\includegraphics[height=1.1in]{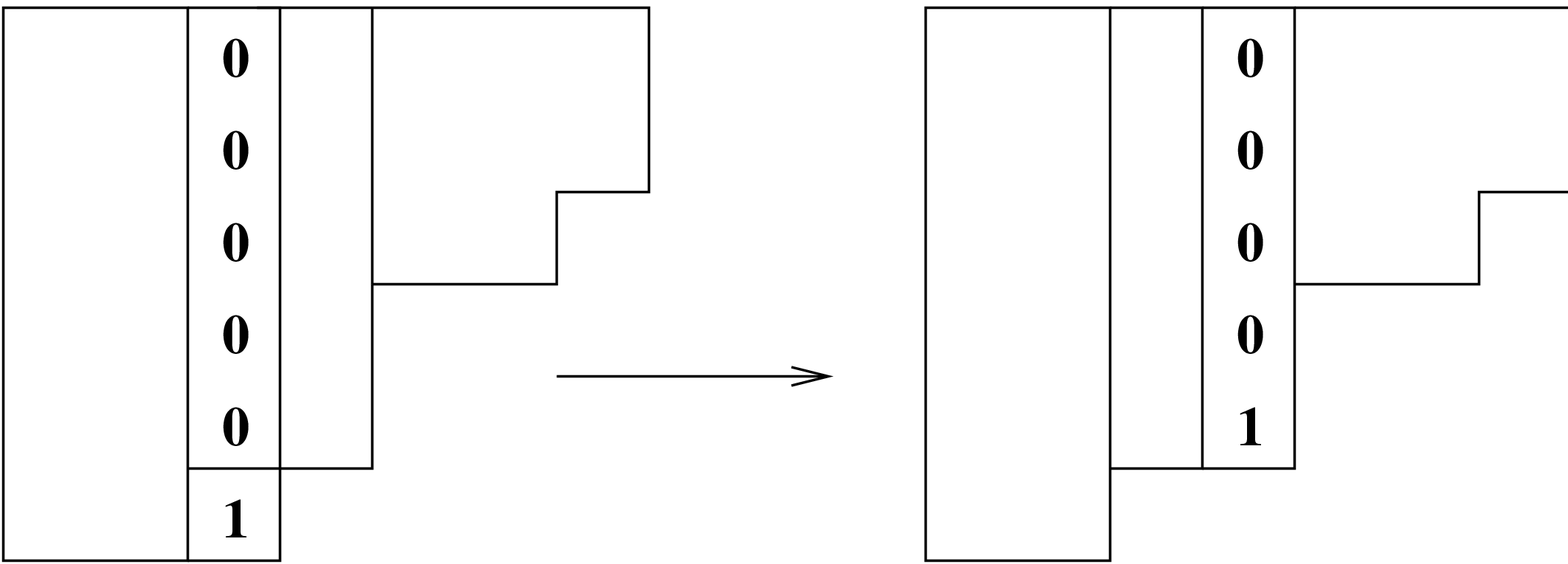}
\caption{}
\label{Right3}
\end{figure}

Suppose that the rightmost entry of the $j$th row of $\SSS$
is a necessary $1$.
Then we define a new tableau $\T$ by
deleting the column containing the necessary $1$ and adding a new
column whose length is $1$ less.  That new column consists entirely
of $0$'s except for a necessary $1$ at the bottom, and it is inserted
as far east a possible.  Such a column cannot introduce a forbidden
pattern (condition (2) in the definition of permutation tableau),
so $\T$ is in fact a permutation tableau.
See Figure \ref{Right3}.

We define $\prob(\SSS \to \T) = \frac{1}{N+1}$.
If the column in $\SSS$ containing the necessary $1$ has length at least $2$
then $\wt(\T) = \wt(\SSS)$, because we have not changed the 
number of superfluous $1$'s, the $1$'s in the top row, or the 
number of unrestricted rows.  Thus  
$\wt(\SSS) \cdot \prob(\SSS \to \T) = \frac{\wt(\T)}{N+1}$.

In the special case that 
the column in $\SSS$ containing the necessary $1$ has length exactly $2$
then $\wt(\T) = {\alpha}^{-1} \wt(\SSS)$, because the new tableau
will have a new $1$ in the top row.  Thus  
$\wt(\SSS) \cdot \prob(\SSS \to \T) = \frac{\alpha \wt(\T)}{N+1}$.

\subsection{Particle exits to the right}\label{exitright}

\begin{figure}[h]
\centering
\includegraphics[height=1.1in]{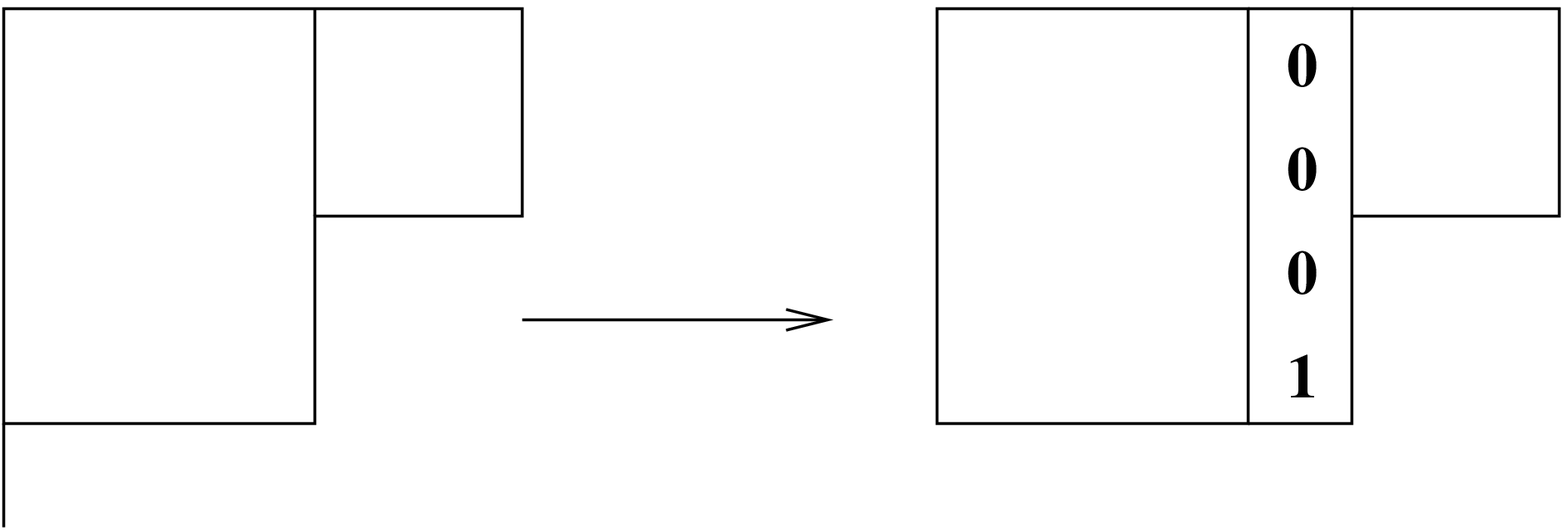}
\caption{}
\label{Exit}
\end{figure}

If $\SSS$ contains a row of length $0$ then 
there is a transition in PT from $\SSS$ that corresponds to
a particle in $\pr(\SSS)$ exiting the PASEP to the right. 
We define a new tableau $\T$ by deleting the $t$th row 
of $\SSS$ (which has length $0$) and adding a new column of 
length $t-1$ which consists of $t-2$ $0$'s followed by a $1$
(read top-to-bottom), inserting this column into the tableau
as far to the right as possible.  The result
is clearly a permutation tableau.  See Figure \ref{Exit}. 

We define $\prob(\SSS \to \T) = \frac{\beta}{N+1}$.
Since we have deleted an unrestricted row (the row of length $0$),
we have $\wt(\T) = \beta \wt(\SSS)$.  Thus  
$\wt(\SSS) \cdot \prob(\SSS \to \T) = \frac{\wt(\T)}{N+1}$.

\subsection{Particle hops left}\label{hopleft}

\begin{figure}[h]
\centering
\includegraphics[height=1.1in]{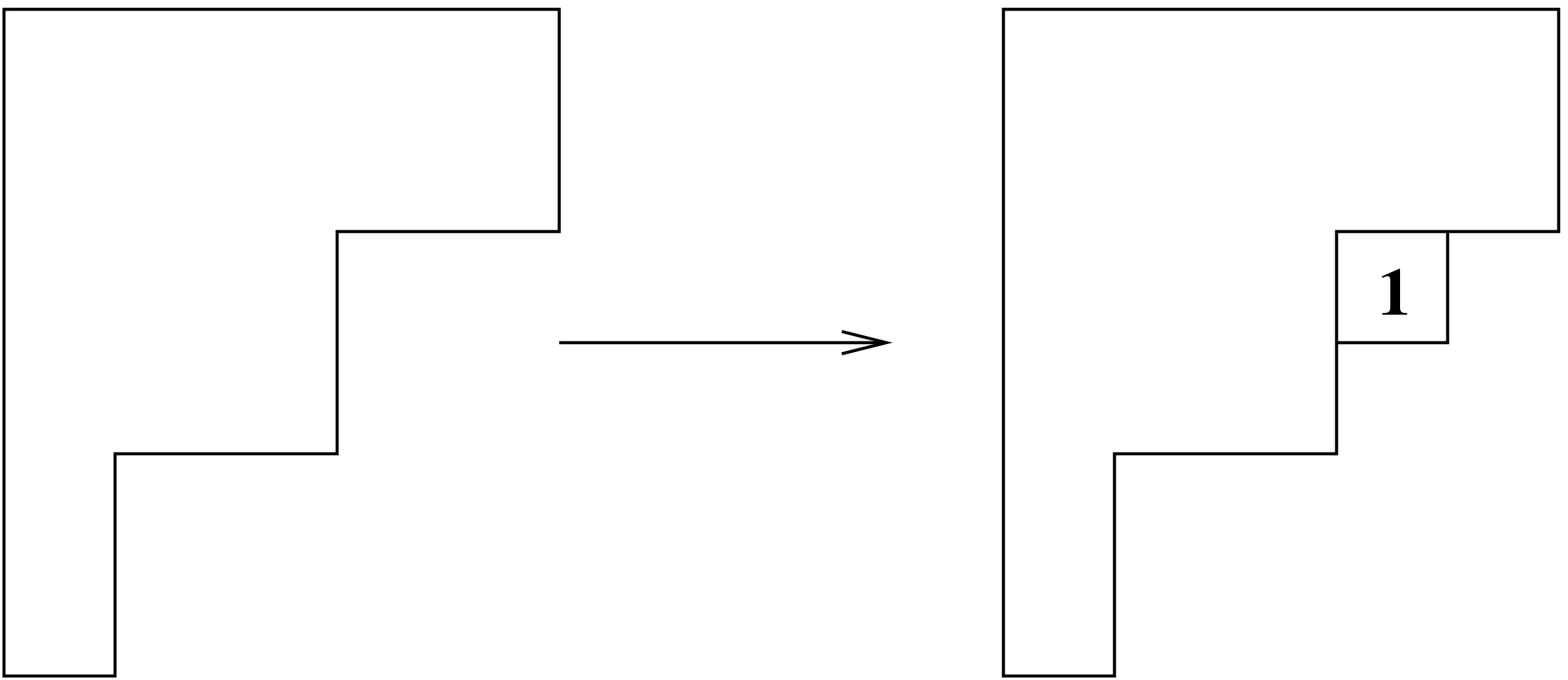}
\caption{}
\label{Left}
\end{figure}

If some row $\lambda_j > \lambda_{j+1}$ in $\SSS$ then there 
is a transition in PT from $\SSS$ that corresponds to the
$j$th black particle in $\pr(\SSS)$ 
hopping to the left in the PASEP.
We define a new tableau $\T$ by increasing the length of 
the $(j+1)$st row by $1$ and filling the extra square with a $1$.
The result is clearly a permutation tableau.
See Figure \ref{Left}.

We define $\prob(\SSS \to \T) = \frac{q}{N+1}$.
Since we have added one superfluous $1$, we have 
that $\wt(\T) = q \wt(\SSS)$.  Thus  
$\wt(\SSS) \cdot \prob(\SSS \to \T) = \frac{\wt(\T)}{N+1}$.

\subsection{Projection of Markov chains}

We now formulate a notion of {\it projection} for Markov chains.  We  
have not been able to find this definition in the literature, but,
for example, the Markov chain of Duchi and Schaeffer \cite{jumping}
is a projection in this sense.

\begin{definition}
Let $M$ and $N$ be Markov chains on finite sets $X$ and $Y$, and 
let $F$ be a surjective map from $X$ to $Y$.  We say that 
$M$ {\rm projects} to $N$ if the following properties hold:
\begin{itemize}
\item If $x_1$ and $x_2$ are in $X$ such that 
$\prob_M(x_1 \to x_2) > 0$, then 
$\prob_M(x_1 \to x_2) = \prob_N(F(x_1) \to F(x_2))$.
\item If $y_1$ and $y_2$ are in $Y$ and $\prob_N(y_1 \to y_2)>0$
then for each $x_1 \in X$ such that $F(x_1)=y_1$ there is a unique
$x_2 \in X$ such that $F(x_2)=y_2$ and 
$\prob_M(x_1 \to x_2) > 0$; moreover, 
$\prob_M(x_1 \to x_2) = prob_N(y_1 \to y_2)$.
\end{itemize}
\end{definition}

Let $\prob_M(x_0 \to x; t)$ denote the probability that 
if we start at state $x_0\in M$
at time $0$, then we are in state $x$ at time $t$.

If $M$ projects to $N$ then a walk on the state diagram
of $M$ is indistinguishable from a walk on the state diagram of $N$ in
the following sense.

\begin{proposition}\label{projectwalk}
Suppose that $M$ projects to $N$.
Let $x_0\in X$ and $y_0, \tilde{y}\in Y$ such that 
$F(x_0)=y_0$.  Then 
\begin{equation*}
\prob_N(y_0 \to \tilde{y}; t) = 
\sum_{\tilde{x} \text{ s.t. }F(\tilde{x})=y} prob_M(x_0\to\tilde{x};t).
\end{equation*}
\end{proposition}

\begin{proof}
We use induction.  The base case $t=0$ is trivially true.  Suppose
that the statement is true for $t-1$.  Let $Y'\subset Y$ be the 
set of all states $y'$ such that there is a transition 
$y'\to \tilde{y}$ with nonzero probability.  Let 
$X' = \{x'\in X \ \vert \ F(x') \in Y'\}$.  By the definition
of projection, these are the only states in $X$ which have a transition
to a state $\tilde{x}$ such that $F(\tilde{x})=\tilde{y}$.  Furthermore,
for each $x_1\in X'$ there exists an $x_2 \in X$ such that 
$F(x_2)=\tilde{y}$ and 
$\prob_M(x_1\to x_2)>0$, and $\prob_M(x_1\to x_2)=
\prob_N(y_1\to y_2)$.

Clearly $\prob_N(y_0\to \tilde{y}; t) = 
\sum_{y'\in Y'} 
\prob_N(y_0\to y'; t-1) \cdot
\prob_N(y'\to \tilde{y})$, 
which by the induction hypothesis and the definition 
of projection is equal to 
\begin{equation*}
\sum_{y'\in Y',} \sum_{x' \text{ s.t. }F(x')=y'} 
\prob_M(x_0\to x'; t-1) \cdot
\prob_M(x'\to \tilde{x}). 
\end{equation*}
This is equal to 
\begin{equation*}
\sum_{x' \text{ s.t. }F(x')\in Y'} 
\prob_M(x_0\to x'; t-1) \cdot
\prob_M(x'\to \tilde{x}). 
\end{equation*}
 But now using again the definition of 
projection, we see that this is equal to $\prob_M(x_0\to x; t)$, 
as desired.
\end{proof}

Proposition \ref{projectwalk} implies the following.

\begin{corollary}\label{projectsteady}
Suppose that $M$ projects to $N$ via the map $F$.   Let $y\in Y$ and 
let $X'=\{x\in X \ \vert \ F(x)=y\}$.  Then the steady state probability
that $N$ is in state $y$ is equal to the steady state 
probabilities that $M$ is in any of the states $x\in X'$.
\end{corollary}

Clearly the operator $\pr$ is a surjective map from the set of 
permutation tableaux of half-perimeter $N+1$ to the states of the PASEP
with $N$ sites.
It is clear from our definition of the PT chain that 
the PT chain projects to the PASEP.

\begin{corollary}
The projection map $\pr$ gives a projection from the PT chain
on permutation tableaux of half-perimeter $N+1$ 
to the PASEP with $N$ sites.
\end{corollary}

\section{Steady state probabilities}\label{Proof}

Our main theorem is the following.

\begin{theorem}\label{main}
Consider the PT chain on permutation tableaux of half-perimeter $N+1$
and fix a permutation tableau $\T$ (of half-perimeter $N+1$).
Then the steady state probability of finding the PT chain 
in state $\T$ is $\frac{\wt(\T)}{\sum_\SSS \wt(\SSS)}$. 
Here, the sum is over all permutation tableaux of half-perimeter 
$N+1$.
\end{theorem}

Since the PT chain projects to the PASEP,
Corollary \ref{projectsteady} and  
Theorem \ref{main} imply Theorem \ref{oldtheorem}.

The goal of this section will be to prove Theorem \ref{main}.  To do so,
we will check the defining recurrences
of the steady state.  
More precisely, it suffices to check the following.  
Fix a state $\T$, let $\Q$ be the collection of all states that 
have transitions to $\T$, and let $\Sset$ be the collection of 
all states that have transitions from $\T$.  Then we need to prove
that 
\begin{equation*}
\wt(\T) = \sum_{Q\in \Q} \wt(Q) \cdot \prob(Q \to \T) + 
          \wt(\T) \cdot (1-\sum_{\SSS\in \Sset} \prob(\T \to \SSS)).
\end{equation*}
The above equation expresses the steady state probability of 
being in state $\T$ in two ways, involving two consecutive times
$t$ and $t+1$; equivalently, it encodes the condition that the transition
matrix has a left eigenvector with eigenvalue $1$.
Combining the terms involving $\wt(\T)$, we get
\def\theequation{E\arabic{equation}}
\begin{equation}\label{mainequation}
\sum_{Q\in \Q} \wt(Q) \cdot \prob(Q \to \T) = 
  \wt(\T) \cdot \sum_{\SSS\in \Sset} \prob(\T \to \SSS).
\end{equation}

In order to check this, 
we will divide the set of states in the PT into four different
classes.  In what follows, let B represent a nonempty string of black
particles
and let W represent a nonempty string
of white particles.  Then we divide the set of permutation tableaux
of half-perimeter $N$ into 
the four classes as follows:

\begin{enumerate}                                                              
\item  $\T$ such that $\pr(\T)$ has the form BWBW...BW.  ($2n$ strings)
\item $\T$ such that $\pr(\T)$ has the form BWBW...BWB. ($2n+1$ strings)
\item $\T$ such that $\pr(\T)$ has the form WBWB...WBW. ($2n+1$ strings)
\item $\T$ such that $\pr(\T)$ has the form WBWB...WB.  ($2n$ strings)
\end{enumerate}

\subsection{Analysis of transitions out of state $\T$}\label{Out}
                                                                                
From a state of type (1), there are $n$ possible transitions of the form
``hop right," and $n-1$ possible transitions of the form ``hop left."
                                                                                
From a state of type (2), there are $n$ possible ``hop right" transitions,
n possible ``hop left" transition, and one ``hop out to the right" transition.
                                                                                
From a state of type (3), there are $n$ ``hop right" transitions, 
$n$ ``hop left"
transitions, and one ``hop in from the left" transition.
                                                                                
From a state of type (4), there are $n-1$ ``hop right" transitions, $n$ 
``hop 
left"
transitions, one ``hop in from the left," and one ``hop out from right."

We can now calculate the quantity
$ \sum_{\SSS\in \Sset} \prob(\T \to \SSS)$ in all four cases.
These quantities are as follows:

\begin{enumerate}
\item $(n+q(n-1))/(N+1)$
\item $(n+qn+\beta)/(N+1)$
\item $(n+qn+\alpha)/(N+1)$
\item $((n-1)+qn+\alpha+\beta)/(N+1)$
\end{enumerate}

\subsection{Analysis of transitions into state $\T$}

\subsubsection{Type (1)}

\begin{figure}[h]
\centering
\includegraphics[height=1.3in]{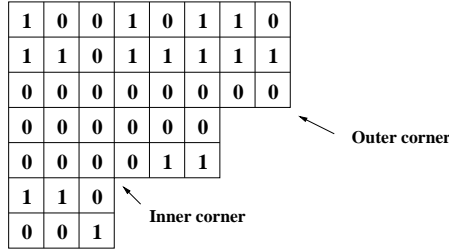}
\caption{A state $\T$ of type (1), where $n=3, a=1,b=1,c=1$}
\label{fig1}
\end{figure}

Fix a state $\T$ of type (1).  The Young diagram of $\T$ will have
$n$ outer corners (south step followed by west step), and $n-1$ inner 
corners (west step followed by south step).  If we say only 
``corner," this will mean an outer corner. 

Of the entries in the $n$ outer corners, $a$ of them are $0$'s, 
$b$ of them are necessary $1$'s, and 
$n-a-b$ are superfluous $1$'s.

For each $0$ corner, there is a state $Q$ such that 
$\prob(Q \to \T)>0$ and such that 
$\wt(Q) \cdot \prob(Q \to \T) = \frac{\wt(\T)}{N+1}$.  Therefore
the $0$ corners contribute $\frac{a \wt(\T)}{N+1}$ to the left-hand-side
of equation \eqref{mainequation}.

For each corner which is a superfluous $1$, 
there is a state $Q$ such that 
$\prob(Q \to \T)>0$ and such that 
$\wt(Q) \cdot \prob(Q \to \T) = \frac{\wt(\T)}{N+1}$.  Therefore
these superfluous $1$'s contribute $\frac{(n-a-b) \wt(\T)}{N+1}$ to the left-hand-side
of equation \eqref{mainequation}.

For each corner which is a necessary $1$, 
there is a state $Q$ such that 
$\prob(Q \to \T)>0$ and such that 
$\wt(Q) \cdot \prob(Q \to \T) = \frac{\wt(\T)}{N+1}$.  Therefore
these necessary $1$'s contribute $\frac{b \wt(\T)}{N+1}$ to the left-hand-side
of equation \eqref{mainequation}.

For each inner corner, 
there is a state $Q$ such that 
$\prob(Q \to \T)>0$ and such that 
$\wt(Q) \cdot \prob(Q \to \T) = \frac{q\wt(\T)}{N+1}$.  Therefore
these inner corners contribute 
$\frac{q(n-1) \wt(\T)}{N+1}$ to the left-hand-side
of equation \eqref{mainequation}.

The sum of all of these contributions is 
$\frac{\wt(\T)}{N+1} (n+q(n-1))$. Comparing this with the results of 
Subsection \ref{Out}, we see that equation \eqref{mainequation}
holds for states of Type (1).

\subsubsection{Type (2)}

\begin{figure}[h]
\centering
\includegraphics[height=1.3in]{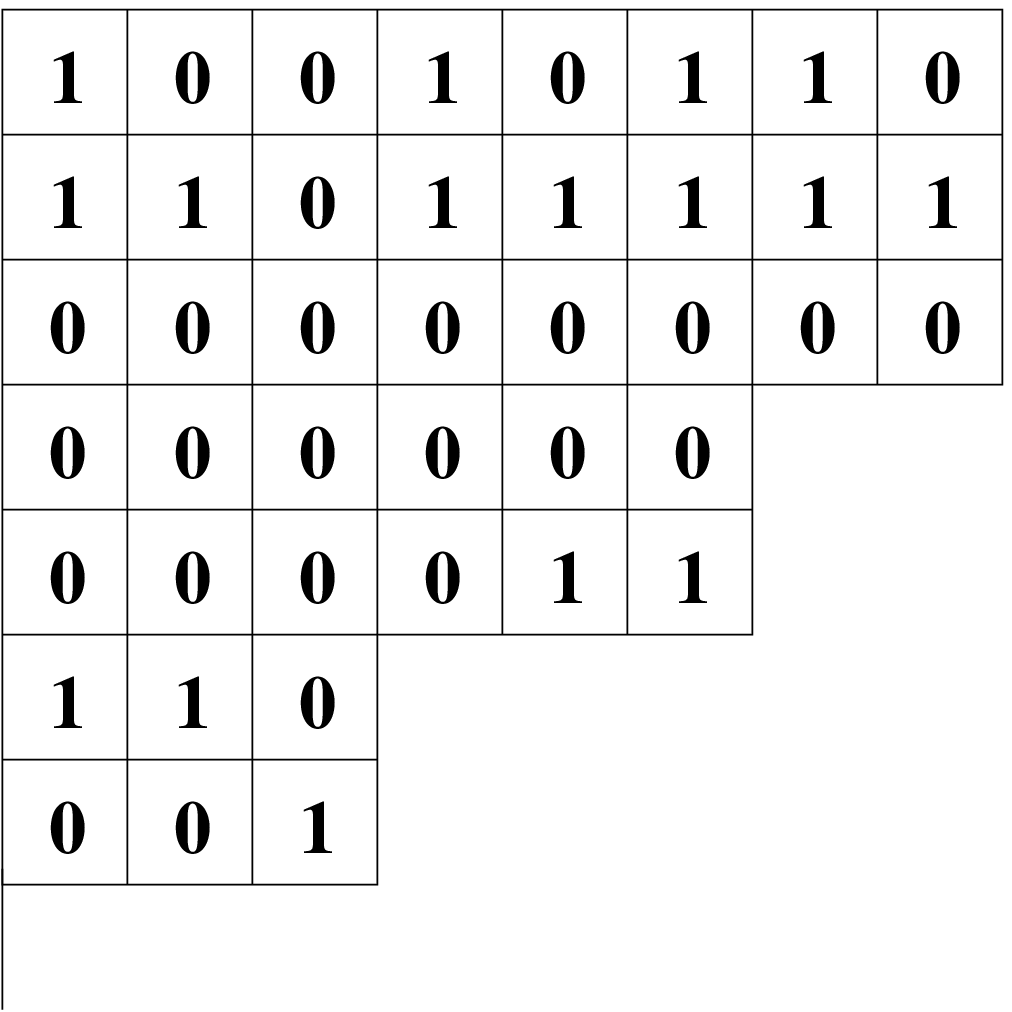}
\caption{A state $\T$ of type (2), where $n=3$}
\label{fig2}
\end{figure}

Fix a state $\T$ of type (2).  The Young diagram of $\T$ will have
$n$ outer corners, and $n$ inner 
corners.

Of the entries in the $n$ outer corners, $a$ of them are $0$'s, 
$b$ of them are necessary $1$'s, and 
$n-a-b$ are superfluous $1$'s.

As before, the $0$ corners 
contribute $\frac{a \wt(\T)}{N+1}$ to the left-hand-side
of equation \eqref{mainequation}.

Similarly, the corners containing
superfluous $1$'s 
contribute $\frac{(n-a-b) \wt(\T)}{N+1}$ to the left-hand-side
of equation \eqref{mainequation}.

The corners containing  necessary $1$'s 
contribute $\frac{b \wt(\T)}{N+1}$ to the left-hand-side
of equation \eqref{mainequation}.

The $n$ inner corners contribute 
$\frac{qn \wt(\T)}{N+1}$ to the left-hand-side
of equation \eqref{mainequation}.

The final empty row of the Young diagram contributes 
$\frac{\beta \wt(\T)}{N+1}$ to equation \eqref{mainequation}.

The sum of all of these contributions is 
$\frac{\wt(\T)}{N+1} (n+qn+\beta)$. Comparing this with what 
we found 
in Subsection \ref{Out}, we see that equation \eqref{mainequation}
holds for states of Type (2).

\subsubsection{Type (3)}

\begin{figure}[h]
\centering
\includegraphics[height=1.3in]{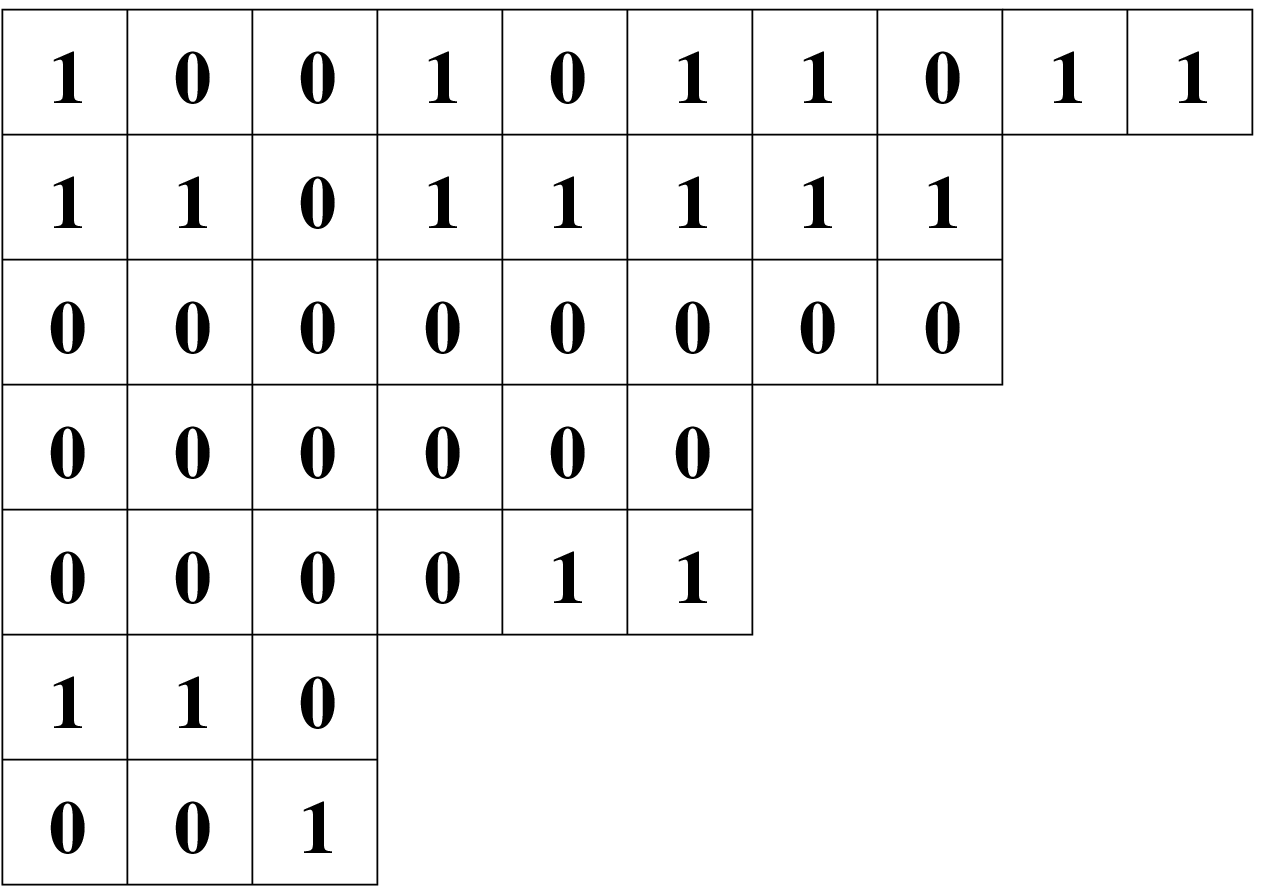}
\caption{A state $\T$ of type (3), where $n=3$}
\label{fig3}
\end{figure}

Fix a state $\T$ of type (3).  The Young diagram of $\T$ will have
$n$ outer corners (not counting the corner formed by the rightmost
column of length $1$), and $n$ inner 
corners.

Of the entries in the $n$ outer corners, $a$ of them are $0$'s, 
$b$ of them are necessary $1$'s, and 
$n-a-b$ are superfluous $1$'s.

As before, the $0$ corners 
contribute $\frac{a \wt(\T)}{N+1}$ to the left-hand-side
of equation \eqref{mainequation}.

The corners containing
superfluous $1$'s 
contribute $\frac{(n-a-b) \wt(\T)}{N+1}$ to the left-hand-side
of equation \eqref{mainequation}.

The corners containing  necessary $1$'s 
contribute $\frac{b \wt(\T)}{N+1}$ to the left-hand-side
of equation \eqref{mainequation}.

The $n$ inner corners contribute 
$\frac{qn \wt(\T)}{N+1}$ to the left-hand-side
of equation \eqref{mainequation}.

The (leftmost) column of length $1$  contributes 
$\frac{\alpha \wt(\T)}{N+1}$ to equation \eqref{mainequation}.

The sum of all of these contributions is 
$\frac{\wt(\T)}{N+1} (n+qn+\alpha)$. Comparing this with  
Subsection \ref{Out}, we see that equation \eqref{mainequation}
holds for states of Type (3).

\subsubsection{Type (4)}

\begin{figure}[h]
\centering
\includegraphics[height=1.3in]{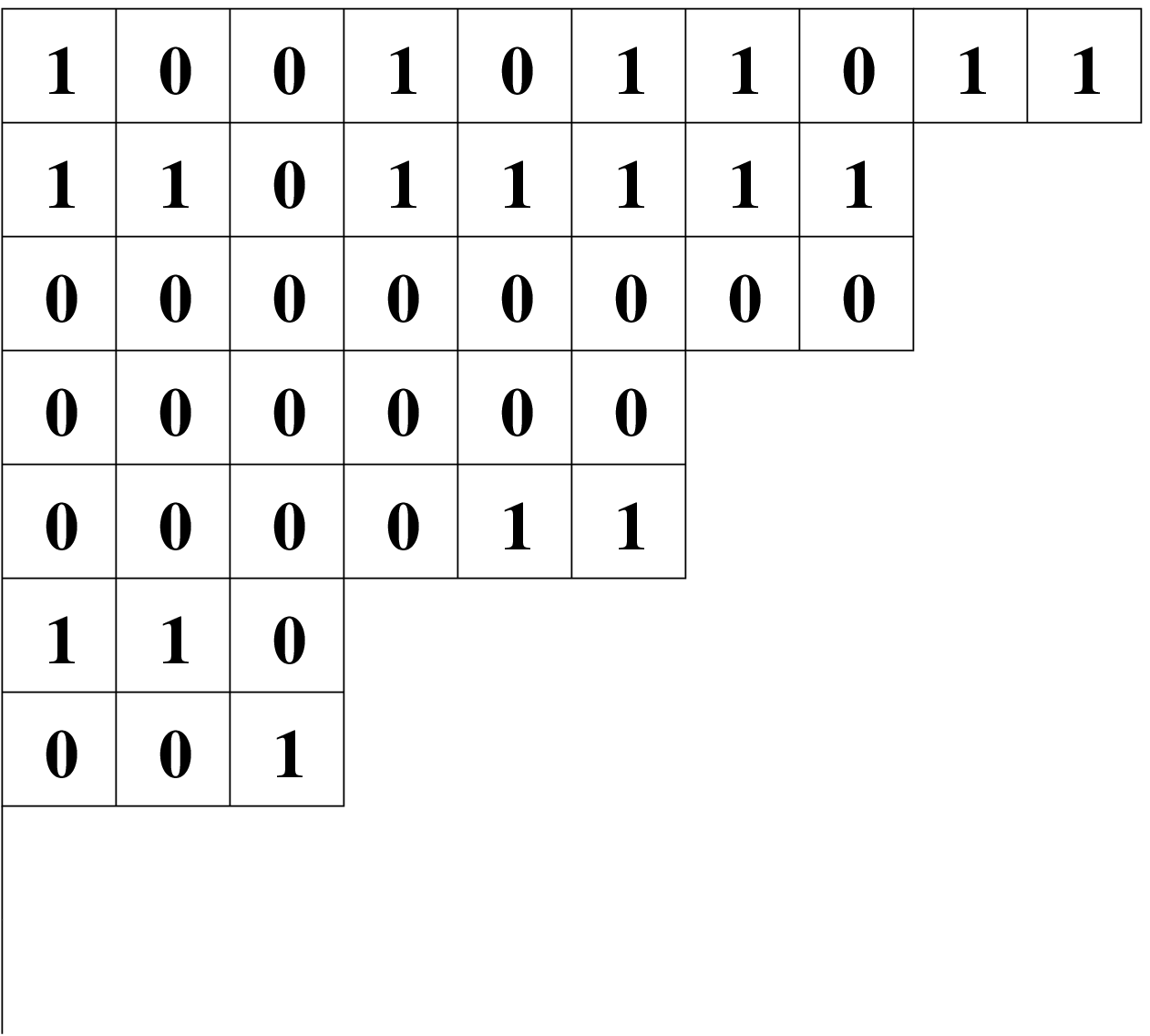}
\caption{A state $\T$ of type (4), where $n=4$}
\label{fig4}
\end{figure}

Fix a state $\T$ of type (4).  The Young diagram of $\T$ will have
$n-1$ outer corners (not counting the corner formed by the rightmost
column of length $1$), and $n$ inner 
corners.

Of the entries in the $n-1$ outer corners, $a$ of them are $0$'s, 
$b$ of them are necessary $1$'s, and 
$n-1-a-b$ are superfluous $1$'s.

As before, the $0$ corners 
contribute $\frac{a \wt(\T)}{N+1}$ to the left-hand-side
of equation \eqref{mainequation}.

The corners containing
superfluous $1$'s 
contribute $\frac{(n-1-a-b) \wt(\T)}{N+1}$ to the left-hand-side
of equation \eqref{mainequation}.

The corners containing necessary  $1$'s 
contribute $\frac{b \wt(\T)}{N+1}$ to the left-hand-side
of equation \eqref{mainequation}.

The $n$ inner corners contribute 
$\frac{qn \wt(\T)}{N+1}$ to the left-hand-side
of equation \eqref{mainequation}.

The (leftmost) column of length $1$  contributes 
$\frac{\alpha \wt(\T)}{N+1}$ to equation \eqref{mainequation}.

The final empty row of the Young diagram contributes 
$\frac{\beta \wt(\T)}{N+1}$ to equation \eqref{mainequation}.

The sum of all of these contributions is 
$\frac{\wt(\T)}{N+1} ((n-1)+qn+\alpha+\beta)$. Comparing this with  
Subsection \ref{Out}, we see that equation \eqref{mainequation}
holds for states of Type (4).

This completes the proof of the main theorem.

\section{The PT chain as a Markov chain on permutations}\label{Perm-Section}

There is a bijection $\Phi$ 
between the set of permutation tableaux of half-perimeter $n$ and 
the permutations in $S_n$ \cite{SW}, which translates various statistics on 
permutation tableaux into statistics on permutations.  This bijection allows us
to interpret the PT chain as a Markov chain on permutations.  

In this section we will describe the PT chain 
in terms of permutations.  First we need to recall the definition of 
the bijection from \cite{SW}.

\subsection{The bijection from permutation tableaux to permutations}

Before defining $\Phi$, it is necessary to introduce some notation.
A {\em weak excedance} of a
permutation $\pi$ is an index $i$  such that $\pi(i) \geq i$.  
A non-excedance is an index $i$ such that $\pi(i) < i$.

We define the {\it diagram}
$D(\T)$ associated with a permutation tableaux $\T$ of shape
$\lambda$ as follows.  
Recall that $p(\lambda)$ is the lattice path which cuts out
the south-east border of $\lambda$.
From north-east to south-west,
label each of the (unit) steps in this path with a number from $1$ to
$n$.
Then, remove the $0$'s from $\T$ and replace each $1$ in $\T$ with a vertex. 
 Finally, from each vertex $v$, draw an
edge to the east and an edge to the south; each such edge should
connect $v$ to either a closest vertex in the same row or column, or
to one of the labels from $1$ to $n$.  The resulting picture is the
{\em diagram} $D(\T)$.  See Figure \ref{fig-diagram}.

We now define the permutation $\pi = \Phi(\T)$ via the following
procedure.  For each $i \in \{1, \dots , n\}$, find the corresponding
position on $D(\T)$ which is labeled by~$i$.  If the label $i$ is on a
vertical step of $P$, start from this position and travel straight
west as far as possible on edges of $D(\T)$. Then, take a ``zig-zag''
path southeast, by traveling on edges of $D(\T)$ south and east and
turning at each opportunity (i.e. at each new vertex).  This path will
terminate at some label $j\ge i $, and we define $\pi(i) = j$.  If $i$ is
not connected to any edge (equivalently, if there are no vertices in
the row of $i$) then we set $\pi(i)=i$.  Similarly, if the label $i$
is on a horizontal step of~$P$, start from this position and travel
north as far as possible on edges of $D(\T)$. Then, as before, take a
zig-zag path south-east, by traveling on edges of $D(\T)$ east and
south, and turning at each opportunity.  This path will terminate at
some label $j<i$, and we let $\pi(i) = j$.

See Figure \ref{zigzag} for a picture of the path taken by $i$.

\psset{unit=.7pt, arrowsize=7pt, linewidth=1pt}
\psset{linecolor=blue}
\newgray{grayish}{.90}
\newrgbcolor{embgreen}{0 .5 0}
\def\vblack(#1, #2)#3{\cnode*[linecolor=black](#1, #2){3}{#3}}
\def\vwhite(#1,#2)#3{\cnode[linecolor=black,fillcolor=white,fillstyle=solid](#1,#2){3}{#3}}
\countdef\x=23
\countdef\y=24
\countdef\z=25
\countdef\t=26

\def\tbox(#1,#2)#3{
\x=#1 \y=#2
\multiply\x by 36
\multiply\y by 36
\z=\x \t=\y
\advance\z by 36
\advance\t by 36
\psline(\x,\y)(\x,\t)(\z,\t)(\z,\y)(\x,\y)
\advance\x by 18
\advance\y by 18
\rput(\x,\y){{\bf #3}}}

\def\evbox(#1,#2)#3{
\x=#1 \y=#2
\multiply\x by 36
\multiply\y by 36
\advance\x by 10
\advance\y by 18
\rput(\x,\y){{\bf #3}}}

\def\ehbox(#1,#2)#3{
\x=#1 \y=#2
\multiply\x by 36
\multiply\y by 36
\advance\x by 18
\advance\y by -12
\rput(\x,\y){{\bf #3}}}

\def\contournumbers{
\evbox(3,0){5}
\evbox(4,1){3}
\evbox(4,2){2}
\evbox(4,3){1}
\ehbox(0,0){8}
\ehbox(1,0){7}
\ehbox(2,0){6}
\ehbox(3,1){4}
}

\def\bdot{\pscircle*{1.5mm}}
\def\wdot{\pscircle{1.5mm}}

\begin{figure}
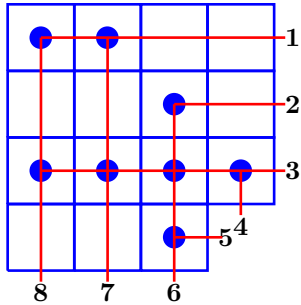

\pspicture(10,-10)(144,160)
\tbox(0,0){}
\tbox(1,0){}
\tbox(2,0){\bdot}

\tbox(0,1){\bdot}
\tbox(1,1){\bdot}
\tbox(2,1){\bdot}
\tbox(3,1){\bdot}

\tbox(0,2){}
\tbox(1,2){}
\tbox(2,2){\bdot}
\tbox(3,2){}

\tbox(0,3){\bdot}
\tbox(1,3){\bdot}
\tbox(2,3){}
\tbox(3,3){}

\contournumbers

\psset{linecolor=red}

\psline(150,126)(18,126)
\psline(18,126)(18,-6)

\psline(150,54)(18,54)
\psline(54,126)(54,-6)

\psline(90,90)(90,-6)
\psline(90,18)(116,18)

\psline(90,90)(150,90)
\psline(126,54)(126,30)

\endpspicture
\caption{\label{fig-diagram} The diagram of a tableau.}
\end{figure}

\begin{figure}
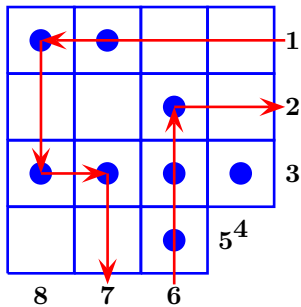

\pspicture(10,-10)(144,160)
\tbox(0,0){}
\tbox(1,0){}
\tbox(2,0){\bdot}

\tbox(0,1){\bdot}
\tbox(1,1){\bdot}
\tbox(2,1){\bdot}
\tbox(3,1){\bdot}

\tbox(0,2){}
\tbox(1,2){}
\tbox(2,2){\bdot}
\tbox(3,2){}

\tbox(0,3){\bdot}
\tbox(1,3){\bdot}
\tbox(2,3){}
\tbox(3,3){}

\contournumbers

\psset{linecolor=red}

\psline{->}(150,126)(18,126)
\psline{->}(18,126)(18,54)
\psline{->}(18,54)(54,54)
\psline{->}(54,54)(54,-6)

\psline{->}(90,-6)(90,90)
\psline{->}(90,90)(150,90)

\endpspicture
\caption{\label{zigzag}The paths taken by 1 and 6: $\pi(1)=7$, $\pi(6)=2$.}
\end{figure}

\begin{example}
If $\T$ is the permutation tableau whose diagram is
given in Figures \ref{fig-diagram} and \ref{zigzag}, then $\Phi(\T) =
74836215$.
\end{example}

Various properties of $\Phi$ were proved in \cite{SW}; we now recall those that will 
be useful to us.

\begin{proposition}\cite{SW}
\begin{itemize}
\item 
In $\Phi(\T)$, the letter $i$ is a fixed point if and only if there is
an entire row in $\T$ that has no 1's and whose right hand edge is
labeled by $i$.  
\item 
The weak excedances  of $\pi = \Phi (\T)$ are precisely the
labels on the vertical edges of $P$.  The non-excedances
of $\pi$ are
precisely the labels on the horizontal edges of $P$.  In particular,
$\Phi(\T)$ is a permutation in $S_n$ with precisely $k$ weak
excedances.
\end{itemize}
\end{proposition}

\subsection{The PT chain on permutations}

Now that we have defined $\Phi$, it is a straightforward exercise to 
translate the PT chain on permutation-tableaux into a Markov chain 
on permutations.  

For convenience, we will first make two definitions.
Consider a permutation $\pi$ on a set $S$.  Let $i$ be in $S$ such that
$\pi(i) \neq i$.  We define the
{\it collapse of $\pi$ at $i$} to be the permutation $\sigma$
on the set 
$S \backslash{\{i\}}$ defined by
$\sigma(\pi^{-1}(i)):=\pi(i)$ and $\sigma(j)=\pi(j)$
for $j \neq \pi^{-1}(i)$. 
For example if $S=\{3,4,6,7\}$ and $\pi=(6,7,4,3)$,
then after collapsing $\pi$ at $7$ we get
$\sigma = (6,3,4)$, a permutation on the set
$\{3,4,6\}$.

We also want the notion of normalizing a permutation.
If $\pi$ is a permutation on an 
ordered set $S=\{s_1<s_2<...<s_{n}\}$ of cardinality
$n$, 
then the {\it normalization of $\pi$} $\Norm(\pi)$ is the permutation 
on $\{1,2,...,n\}$ that we
get by replacing $s_i$ with $i$.

\begin{figure}[h]
\centering
\includegraphics[height=3.5in]{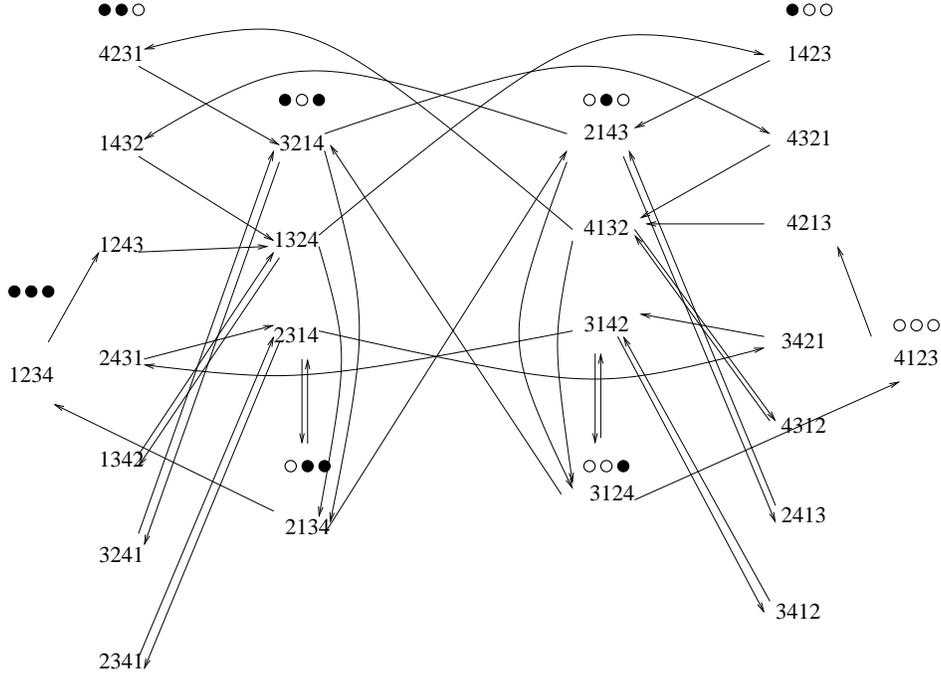}
\caption{The state diagram of the PT chain, depicted in terms of permutations}
\label{states2}
\end{figure}

We are now ready to describe the PT chain in terms of permutations.
An example is shown in Figure \ref{states2}, which is 
simply the result of applying $\Phi$ to Figure \ref{Markov-states}.

In what follows, we will use the notation 
$\{a, \widehat{b}, c\}$ to denote the set $\{a,c\}$.

\begin{proposition}
After applying $\Phi$, the PT chain is the Markov chain on 
the permutations $S_{N+1}$ with the following transitions.
In all cases, $\pi \in S_{N+1}$.


$\bullet$ ``Particle enters from the left"

Suppose that $\pi(2)=1$.  Let $i+1$ be the minimal number greater than
$2$ which is a non-excedance in $\pi$.   
Define $\sigma$ to be the permutation on 
$\{1, \widehat{2}, 3, \dots , N+1\} \cup \{i+\frac{1}{2}\}$
obtained by collapsing $\pi$ at $2$, and adding a fixed point at 
$i+\frac{1}{2}$.
In this case there is a transition $\pi \to \Norm(\sigma)$ 
with probability 
$\prob(\pi\to \Norm(\sigma))=\frac{\alpha}{N+1}$.

$\bullet$ ``Particle hops right"

Suppose that for some $i \geq 2$, $\pi(i) \geq i$ 
and $\pi(i+1)<i+1$.
Define $\sigma$ as follows:

Case 1: If $\pi(i) = i$ then
let $j+1$ be the smallest non-excedance  
of $\pi$ such that $j>i$.  Define $\sigma$ to be the permutation on 
$\{1, 2, \dots , i-1,\widehat{i}, i+1, \dots, N+1\} \cup \{i+\frac{1}{2}\}$
obtained by collapsing $\pi$ at $i$ and inserting a 
new fixed point at $(j+1/2)$.

Case 2: If $\pi(i)>i$ and 
$\pi(i+1)<i$, then define $\sigma$ by 
$\sigma(i)=\pi(i+1)$, $\sigma(i+1)=\pi(i)$,
and $\sigma(j)=\pi(j)$ for $j \neq i, i+1$.

Case 3: If $\pi(i)>i$ and 
$\pi(i+1)=i$, then 
let $j$ be the greatest
number less than $i$ which is a weak excedance 
(it exists since $1$ is a weak excedance).  Let
$b = \sigma^{-1}(j)$.  Let $\sigma$ be the permutation
on $\{1, 2, \dots , \widehat{i+1}, \dots, N+1 \} \cup \{j+\frac{1}{2}\}$
obtained by collapsing $\pi$ at $i+1$, and replacing 
$\pi(b)=i$ by $\sigma(b)=i+\frac{1}{2}$ and $\sigma(i+\frac{1}{2})=i$.

Then there is a transition $\pi \to \Norm(\sigma)$ with probability
$\frac{1}{N+1}$.

$\bullet$ ``Particle exits to the right"

If $\pi(N+1)=N+1$, then let 
$i$
be the maximal number less than $N+1$ such that
$\pi(i) \geq i$, and let 
$a=\pi^{-1}(i)$.  Define $\sigma$ to be the permutation on 
$\{1, 2, \dots , N\} \cup \{i+\frac{1}{2}\}$ obtained by collapsing
$\pi$ at $N+1$ and replacing $\pi(a)=i$ by 
$\sigma(a)=i+\frac{1}{2}$ and $\sigma(i+\frac{1}{2})=i$.
Then there is a transition $\pi \to \Norm(\sigma)$ such that
$\prob(\pi\to\Norm(\sigma))=\frac{\beta}{N+1}$.

$\bullet$ ``Particle hops left"

If $\pi(i)<i$ and $\pi(i+1) \geq i+1$ then define
$\sigma$ to be the permutation on $\{1, 2, \dots , N+1\}$
defined by 
$\sigma(i)=\pi(i+1)$, $\sigma(i+1)=\pi(i)$ and $\sigma(j)=\pi(j)$
for $j \neq i, i+1$.  
Then there is a transition $\pi \to \sigma$ such that 
$\prob(\pi\to\sigma)=\frac{q}{N+1}$.
\end{proposition}

\begin{proof}
This proof follows easily from the definition of the PT chain and the 
bijection $\Phi$.  Clearly inserting a new all-zero row into a tableau
corresponds to inserting a new fixed point into a permutation, i.e.\
inserting a ``minimal weak excedance" $i \to i$.  
Inserting a new column which consists from top to bottom of zeros
and a single one corresponds to adding a new ``minimal non-excedance"
$i+1 \to i$.  Finally, adding or removing an outer corner to a tableau
corresponds to switching the values of $\pi(i)$ and $\pi(i+1)$.
\end{proof}

Now that we have defined the PT chain in terms of permutations, we
need to define the maps and statistics on permutations that 
are relevant to the PASEP.

We define a surjective map $\pr$ from $S_{N+1}$ to states of the PASEP with
$N$ sites as follows:
if $\pi \in S_{N+1}$ has $W$ as its set of weak-excedances,
then $\pr(\pi) = (a_1, \dots , a_{N}) \in \{0,1\}^N$ is defined by 
$a_i = 1$ if $i+1 \in W$ and $a_i = 0$ otherwise.  Clearly the
map $\pr$ on permutation tableaux corresponds via $\Phi$ to this
map $\pr$ on permutations.

We define a {\it crossing} of a permutation $\pi$ to be a pair of 
indices $i$ and $j$ such that either $i<j \leq \pi(i) < \pi(j)$ or 
$\pi(i)<\pi(j)<i<j$.  In \cite{SW} it was shown that 
the superfluous ones in a permutation tableaux $\T$ are in bijection
with the crossings of $\Phi(\T)$.  Therefore we define
$\rk(\pi)$ to be the number of crossings of $\pi$.

We define a {\it LR-maximum} of $\pi$ to be
an index $i$ such that $\pi(i)>\pi(j)$ for $j<i$, and we call it 
{\it special} if in addition, $\pi(i)>\pi(1)$.
Similarly we define a {\it RL-minimum} 
of $\pi$ to be an index $i$ such that $\pi(i)< \pi(j)$ for $j>i$,
and we call it {\it special} if in addition, $\pi(i)<\pi(1)$.
Abusing notation, we let
$f(\pi)$ denote the number of special 
RL-minima, and we let 
$u(\pi)$ denote the number of special LR-maxima of $\pi$.
We show in \cite{csw} that 
if $\T$ is a permutation tableau, then 
$f(\Phi(\T)) = f(\T)$ and $u(\Phi(\T)) = u(\T)$.

In analogy with our weight function on tableaux, we now define
$\wt^{q,\alpha,\beta}({\pi})=q^{rk({\pi})}\alpha^{-f({\pi})}
\beta^{-u({\mathcal T})}$.  

We can now translate Theorem \ref{main} into the language of
permutations.

\begin{corollary}\label{maincorollary}
Consider the PT chain on permutations in  $S_{N+1}$
and fix a permutation $\pi\in S_{N+1}$.
Then the steady state probability of finding the PT chain 
in state $\pi$ is $\frac{\wt(\pi)}{\sum_{\pi'} \wt(\pi')}$. 
Here, the sum is over all permutations $\pi'\in S_{N+1}$. 
\end{corollary}

\section{The involution}\label{involution}

The goal of this section is to introduce an involution $I$ on
permutations (equivalently, on permutation tableaux) which generalizes the 
particle-hole symmetry of the PASEP, and reveals a symmetry in the PT chain.  
This symmetry 
is a graph automorphism of the state diagram -- this is 
depicted as a reflective symmetry from left to right 
in Figure \ref{Markov-states}.
Using the same notation that we did for the particle-hole symmetry on states of the PASEP, 
we will use $I(\pi)$ or $\overline{\pi}$, or 
$I(\T)$ or $\overline{\T}$ to denote the image of $\pi$ or $\T$ under
the involution.
Because in this case it is 
easier to work with permutations, we will give all proofs in terms
of permutations, and only {\it define} the analogous involution
on permutation tableaux. Philippe Duchon has informed us that he independently
discovered such an involution \cite{du}.

\begin{definition}
Let $\pi=(\pi(1),\ldots ,\pi(N+1))$ be an element of 
$S_{N+1}$.  We define the permutation 
$\overline{\pi}=(\overline{\pi}(1),\ldots,\overline{\pi}(N+1))$
as follows:
\begin{align*}
\overline{\pi}(1)&=N+2-\pi(1)\\  
\overline{\pi}(i)&=N+2-\pi(N+3-i), \ \ \ 2\le i\le N+1.
\end{align*}
\end{definition}
We will prove the following result.

\begin{theorem}\label{involution-theorem}
The map $I$ on permutations in $S_{N+1}$ which sends
 $\pi$ to $\overline{\pi}$ 
has the following properties:
\begin{enumerate}
\item $\pr(\pi) = \overline{\pr(\overline{\pi})}$.  In other words, 
$\pr(\pi)$ and 
   $\pr(\overline{\pi})$ are related via the particle-hole symmetry.
\item  $\rk({\pi})=\rk(\overline{\pi})$.
 \item $u(\pi)=f(\overline{\pi})$.
\item  $f(\pi)=u(\overline{\pi})$.
\end{enumerate}
\end{theorem}

Theorem \ref{involution-theorem} and Corollary \ref{maincorollary} then 
immediately imply the following result. 

\begin{corollary}\label{PH-extension}
Conider the PT chain on permutations in $S_{N+1}$
and fix $\pi \in S_{N+1}$.
Then 
$\wt({\pi},q,\alpha,\beta)=\wt(\overline{{\pi}},q,\beta,\alpha).$
Moreover, the steady state probability of finding the PT chain in state
$\pi$ is equal to the steady state probability of finding the PT chain in state $\overline{\pi}$.
\end{corollary}

This can be seen as an extension of the particle-hole symmetry that was mentioned in Remark 
\ref{particle-hole}.
Indeed our Theorem 3.1 states
that the probability to be in state $\tau$ is the (normalized) weight-generating function
$\frac{F_{\lambda(\tau)}(q)}{Z_N}$
for all permutation tableaux of shape $\lambda(\tau)$.  
After translating this into the corresponding statement for 
permutations, Corollary \ref{PH-extension} immediately
implies that $F_{\lambda(\tau)}(q) = F_{\lambda(\overline{\tau})}(q)$.

Additionally, the following result reveals a symmetry in the state diagram of the PT chain.

\begin{theorem}\label{symmetry}
There is a transition in the PT chain from 
${\pi}$ to ${\sigma}$ if and
only if there is a transition from ${\overline {\pi}}$ to ${\overline {\sigma}}$.  Furthermore, the transition probabilities are related 
as follows:
\begin{enumerate}
\item 
$\prob(\pi \to {\sigma})=\frac{\alpha}{N+1}$ if and only if 
$\prob(\overline{\pi} \to \overline{\sigma})=
\frac{\beta}{N+1}$. 
\item 
$\prob(\pi \to {\sigma})=\frac{\beta}{N+1}$ if and only if 
$\prob(\overline{\pi} \to \overline{\sigma})=
\frac{\alpha}{N+1}$. 
\item 
$\prob(\pi \to {\sigma})=\frac{1}{N+1}$ if and only if 
$\prob(\overline{\pi} \to \overline{\sigma})=
\frac{1}{N+1}$. 
\item 
$\prob(\pi \to {\sigma})=\frac{q}{N+1}$ if and only if 
$\prob(\overline{\pi} \to \overline{\sigma})=
\frac{q}{N+1}$. 
\end{enumerate}
\end{theorem}

\subsection{Proof of Theorem \ref{involution-theorem}}

We present a series of Lemmas that give a refined version of
Theorem \ref{involution-theorem}.

\begin{lemma}\label{involution-lemma}
$I$ is an involution.
\end{lemma}
\begin{proof}
We have $\overline{\overline{\pi}}(1) = N+2 - (N+2-\pi(1)) = \pi(1)$, and 
for $i \geq 2$, $\overline{\overline{\pi}}(i) = 
  N+2 - \overline{\pi}(N+3-i) = N+2-(N+2-\pi(N+3-(N+3-i))) = 
\pi(i)$.
\end{proof}

\begin{lemma} 
Let $\pi\in S_{N+1}$.  Then for $i \geq 2$, $i$ is a weak excedance
of $\pi$ if and only if $N+3-i$ is not a weak excedance of $\overline{\pi}$.
It follows that
$\pr(\pi) = \overline{\pr(\overline{\pi})}$. 
\label{lem0}
\end{lemma}

\begin{proof}  
By definition, $\overline{\pi}(N+3-i) = N+2-\pi(i)$.  Note that
$\pi(i) \geq i$ if and only if 
$\overline{\pi}(N+3-i) \leq N+2-i$, i.e.\ $N+3-i$ is {\it not}
a weak excedance of $\overline{\pi}$.  Since the projection operator
from permutations to states of the PASEP forgets about $\pi(1)$ and 
maps the permutation to a state based on whether the other 
$\pi(i)$ are weak excedances, what we have proved implies that 
$\pr(\pi) = \overline{\pr(\overline{\pi})}$. 
\end{proof}

\begin{lemma}
For $i$ and $j$ not equal to $1$, there is a crossing in positions
$i$ and $j$ in $\pi$ if and only if there is a crossing in positions
$N+3-i$ and $N+3-j$ in $\overline{\pi}$.  
Additionally, the number of crossings involving position $1$ in $\pi$
is equal to the number of crossings involving position $1$ in $\overline{\pi}$.
Therefore $\rk(\T) = \rk(\overline{\T})$.
\label{lem1}
\end{lemma}

\begin{proof}
Consider a crossing in positions $i$ and $j$ in $\pi$ such that 
neither $i$ nor $j$ is $1$.  Without loss of generality, assume
that $i<j\leq \pi(i) < \pi(j)$.  Then it follows that 
$N+2-\pi(j)<N+2-\pi(i)<N+3-j<N+3-i$.  But since 
$\overline{\pi}(N+3-i)=N+2-\pi(i)$ and
$\overline{\pi}(N+3-j)=N+2-\pi(j)$, we have that 
$\overline{\pi}(N+3-j)<\overline{\pi}(N+3-i)<N+3-j<N+3-i$ which 
implies that $N+3-i$ and $N+3-j$ are positions of a crossing in 
$\overline{\pi}$.

The number of crossings involving
position $1$ in ${\pi}$ is equal
to the number of indices $k$ such that 
$1 < k \leq {\pi}(1)<\pi(k)$. Therefore
the number of  crossings involving
position $1$ in $\overline{\pi}$ is equal to
\begin{eqnarray*}
\# \{k \ \vert \ k \leq \overline{\pi}(1) < \overline{\pi}(k) \}&=&
\# \{k \ \vert \ k \leq N+2-{\pi(1)} < N+2-{\pi}(N+3-k) \}\\
&=& \# \{k \ \vert \ k>{\pi(1)}>{\pi}(k) \}\\
&=& \# \{k \ \vert \ k\le {\pi(1)}<{\pi}(k) \}.\\
\end{eqnarray*}
The first step comes from the definition of $\overline\pi$. The second step
comes from replacing $k$ by $N+3-k$, decreasing by $N+2$ and negating. 
The last step comes from the fact that $\pi$ is a permutation.
\end{proof}

\begin{lemma}\label{lem2}
The index $i$ is a special RL-minimum of $\pi$ if and only if 
the index $N+3-i$ is a special LR-maximum of $\overline{\pi}$.
And the index $i$ is a special LR-maximum of $\pi$ if and only if 
the index $N+3-i$ is a special RL-minimum of $\overline{\pi}$.
\end{lemma}

\begin{proof}
We just prove the first part. The second part follows thanks to the
Lemma \ref{involution-lemma}.
Suppose that the index $i$ is a special RL-minimum of $\pi$.  This 
means that $\pi(i)<\pi(1)$
and $\pi(i) < \pi(j)$ for $j>i$.  By definition, 
$\overline{\pi}(1) = N+2-\pi(1)$ and 
$\overline{\pi}(N+3-i) = N+2-\pi(N+3-(N+3-i)) = N+2 - \pi(i)$.
So $\pi(i) < \pi(1)$ implies that 
$\overline{\pi}(N+3-i) = N+2-\pi(i) > N+2-\pi(1)=\overline{\pi}(1)$.
Now note that since $\pi(i)<\pi(j)$ for $j>i$ we have
$\overline{\pi}(N+3-i) =N+2-\pi(i)>N+2-\pi(j) = \overline{\pi}(N+3-j)$,
which completes the proof that $N+3-i$ is a special 
LR-maximum of 
$\overline{\pi}$.
\end{proof}

The combination of Lemmas \ref{involution-lemma}-\ref{lem2} 
therefore implies Theorem \ref{involution-theorem}.

\subsection{Symmetry in the PT chain}
\label{MC}
In this section we will prove Theorem \ref{symmetry}.

\begin{proof}[Proof of Theorem \ref{symmetry}]
Let $\pi$ be a permutation in $S_{N+1}$. 
We will analyze in turn 
the various kinds of transitions from ${\pi}$ to ${\sigma}$ in the 
PT chain.

\begin{itemize}
\item ``Particle enters from the left" and ``exits to the right": 

Note that $\pi(2)=1$ if and only if $\overline{\pi}(N+1)=N+1$,
so there is a transition of type ``particle enters from the left"
out of $\pi$ if and only if there is a transition of type 
``particle exits to the right" out of $\overline{\pi}$.

\item ``Particle hops right, Cases 1 and 3":

If $\pi(i)=i$ and $\pi(i+1)<i+1$ then 
$\overline{\pi}(N+2-i)=N+2-\pi(i+1)>N+1-i$
and $\overline{\pi}(N+3-i)=N+2-\pi(i)=
N+2-i$.
Conversely, $\pi(i)>i$ and $\pi(i+1)=i$ then 
$\overline{\pi}(N+2-i)=N+2-\pi(i+1)=N+2-i$ 
and $\overline{\pi}(N+3-i)=N+2-\pi(i)<N+2-i$
which implies that $\overline{\pi}(N+3-i)<N+3-i$.
So there is a Case 1 transition out of $\pi$ 
at index $i$ if and only if there is a 
Case 3 transition out of $\overline{\pi}$
at $N+2-i$, and vice-versa.

\item  ``Particle hops right, Case 2":

If $\pi(i)>i$ and $\pi(i+1)<i$ then 
$\overline{\pi}(N+3-i)=N+2-\pi(i)<N+2-i$ and 
$\overline{\pi}(N+2-i)=N+2-\pi(i+1)>N+2-i$.  So there is a 
``particle hops right Case 2" transition out of $\pi$ if and only if
there is a ``particle hops right Case 2" transition out of 
$\overline{\pi}$.

\item ``Particle hops left": 

If $\pi(i)<i$ and $\pi(i+1) \geq i+1$ then 
$\overline{\pi}(N+3-i)=N+2-\pi(i)>N+2-i$ implies that
$\overline{\pi}(N+3-i)\geq N+3-i$ and 
$\overline{\pi}(N+2-i)=N+2-\pi(i+1) \geq N+1-i$ implies
that $\overline{\pi}(N+2-i) < N+2-i$.  So there
is a ``particle hops left" transition out of $\pi$ at 
index $i$ if and only if there is a ``particle hops left" transition
out of $\overline{\pi}$ at index $N+2-i$.

\end{itemize}

We also need to show that the 
transition probabilities
$\prob(\pi \to \sigma)$
and $\prob(\overline{\pi} \to
\overline{\sigma})$ are 
related as specified in 
Theorem \ref{symmetry}.
These calculations are straightforward
but tedious so we will show only
the first case.

Suppose $\pi(2)=1$ and 
let $\sigma$ be the permutation obtained after
a transition of the form 
``particle enters from the left." 
Let $i+1$ be the minimum number greater
than $2$ which is a nonexcedance.
So $\sigma$ is the permutation
on $\{1,\hat{2},3,\dots,N+1\}
\cup \{i+\frac{1}{2}\}$ obtained
by collapsing $\pi$ at $2$
and adding a fixed point at 
$i+\frac{1}{2}$.  

Let $\mu$ be the permutation
obtained from $\overline{\pi}$
after a transition of the form
``particle exits to the right."
Let $j$ be the maximum number
less than $N+1$ such that 
$\overline{\pi}(j) \geq j$.
So $j=N+2-i$.
Let $a=\overline{\pi}^{-1}(N+2-i)$.
Then $a=N+3-\pi^{-1}(i)$.
Then $\mu$
is the permutation on 
$\{1,2,\dots,N\} \cup \{N+2\frac{1}{2}-i\}$
obtained by collapsing at $N+1$ and replacing
$\overline{\pi}(a)=N+2-i$ by 
$\mu(a)=N+2\frac{1}{2}-i$
and $\mu(N+2\frac{1}{2}-i)=N+2-i$.
It is now straightforward to see that in fact
$\mu = \overline{\sigma}.$

\end{proof}

We now define the involution in terms of permutation tableaux.

Let ${\T}$ be a permutation tableau with $K$ rows and $N+1-K$ columns
and shape $\lambda=(\lambda_1, \dots, \lambda_K)$.
Numbering the rows from top to bottom and the columns from left to right,
let ${\T}(i,j)$ denote the filling of the cell $(i,j)$ of ${\mathcal T}$.
The {\it conjugate} of ${\mathcal T}$, which we shall denote by 
${\mathcal T}'$, is the tableau of shape $\lambda'$ such that ${{\mathcal T}'}(i,j)={\mathcal T}(j,i)$
for all $i,j$.  Here $\lambda'=(\lambda'_1, \dots , \lambda'_{N+1-K})$ is the
conjugate partition, i.e.\ the partition
formed by the {\it columns} of $\lambda$.
If $a\in \{0,1\}$, let $a^c$ denote $1-a$.

Let $\T$ be a permutation tableau with shape $\lambda = (\lambda_1, \dots , \lambda_K)$
and conjugate shape $\lambda'=(\lambda'_1, \dots , \lambda'_{N+1-K})$.
We define $\overline{\T}$ to be the  tableau of shape
$(K-1, \lambda'_1 - 1, \lambda'_2 - 1, \dots , \lambda'_{N+1-K}-1)$ 
whose entries
are as follows.
\begin{enumerate}
\item $\overline{\mathcal T}(1,j)=1$ if row $j+1$ of ${\mathcal T}$ is unrestricted and 0
otherwise for $1\le j\le K-1$.
\item
$\overline{\mathcal T}(i,j)={{\mathcal T}(j+1,i-1)}^c$ if cell $(j+1,i-1)$ of ${\mathcal T}$
contains a topmost
one or a rightmost restricted zero, and ${\mathcal T}(j+1,i-1)$ otherwise.\\
\end{enumerate}
One can check that
$\pi=\Phi(\T)$ if and only if $\bar\pi=\Phi(\bar\T)$.

\end{document}